\documentclass{article}
\usepackage{graphicx} 
\usepackage{geometry}
\usepackage{amsmath}
\usepackage[utf8]{inputenc}
\usepackage{breqn}
\usepackage{amssymb}
\usepackage{setspace} 
\usepackage{amsfonts}
\usepackage{amsthm}
\usepackage{mathrsfs}
\usepackage{subcaption} 
\usepackage[utf8]{inputenc}
\usepackage{textcomp}
\usepackage{titlesec}
\usepackage{xcolor}
\usepackage{placeins}
\usepackage{tocloft}
\usepackage{comment}

\usepackage{url}
\usepackage[nocompress]{cite}

\usepackage[colorlinks=true, citecolor=violet]{hyperref}

\DeclareMathOperator*{\argmin}{arg \, min}

\newcommand{\R}{\mathbb{R}}
\newcommand{\X}{\mathbf{X}}
\newcommand{\rr}{\mathbf{R}}
\newcommand{\bu}{{\mathbf{u}}}

\newcommand{\uh}{{\hat{u}}}

\newcommand{\revision}[1]{#1}
\newcommand{\revisiontwo}[1]{#1}
\newcommand{\revisionthree}[1]{#1}

\title{Attractor learning for spatiotemporally chaotic dynamical systems using echo state networks with transfer learning}
\author{
Mohammad Shah Alam\thanks{Dept. of Natural and Behavioral Sciences, Sul Ross State University, Eagle Pass, TX 78852, msa24gn@sulross.edu} \quad
William Ott\thanks{Dept. of Mathematics, University of Houston, Houston, TX 77204, 
william.ott.math@gmail.com} \quad
Ilya Timofeyev\thanks{Dept. of Mathematics, University of Houston, Houston, TX 77204, itimofey@Central.uh.edu}
}
\date{\today}

\begin{document}

\maketitle

\begin{abstract} 
\noindent
In this paper, we explore the predictive capabilities of echo state networks (ESNs) for the generalized Kuramoto-Sivashinsky (gKS) equation, an archetypal nonlinear PDE that exhibits spatiotemporal chaos. 
Our research focuses on predicting changes in long-term statistical patterns of the gKS model that result from varying the dispersion relation or the length of the spatial domain.
We use transfer learning to adapt ESNs to different parameter settings and successfully capture changes in the underlying chaotic attractor.
\revisiontwo{Previous work has shown that transfer learning can be used effectively with ESNs for single-orbit prediction.
The novelty of our paper lies in our use of this pairing to predict the long-term statistical properties of spatiotemporally chaotic PDEs.}
\revisiontwo{Nevertheless, we also} show that transfer learning nontrivially improves the length of time that predictions of individual gKS trajectories remain accurate.
\end{abstract}



\section{Introduction}
Over the past decade, there has been a substantial increase in research focused on the application of machine learning to dynamical systems.
Neural networks have played a particularly important role, as a variety of promising techniques based on neural network architectures have been developed.
In particular, reservoir computing~\cite{jaeger2004harnessing, esp1, rc2023survey} has been successfully utilized to learn and predict the dynamics of complex systems, including chaotic models.
Reservoir computing refers to a broad class of recurrent neural networks, wherein the reservoir has its own nonlinear dynamics and only a linear output layer is trained to match the observational data.

In this paper, we explore the predictive capabilities of echo state networks (ESNs) in the context of spatiotemporal chaos.
An echo state network is a simple reservoir computing architecture with the echo state property~\cite{jaeger2001echo,jaeger2002NIPS,maass2002real,buehner2006,jaeger2007echo,esp3, esp1}.
We \revisiontwo{utilize a} methodology that integrates ESNs with transfer learning (TL)~\cite{inubushi2020transfer, chen2022predicting}.

\revision{It is often useful to model a physical process using a nonstationary dynamical system.
Here, the dynamical model itself varies in time and may have time-dependent parameters.
It is challenging to develop machine learning (ML) models that can track parametric changes in nonstationary dynamical systems.
In particular, data in a new parameter regime
might be missing, or insufficient to properly train a de novo ML model.
There exists substantial recent work on developing ML models that can track such time-dependent parametric changes~\cite{patel2021nonstationary, patel2023tipping, pershin2023tipping, bury2021tipping, kong2021tipping, xiao2021bifurcation}.}

\revision{In this paper, we integrate ESNs with transfer learning to predict changes in long-term statistical properties of spatiotemporally chaotic PDEs that result from instantaneously jumping from one parameter regime to another.
Our work here complements existing work on developing ML models that can track parametric changes, as this existing work assumes parameters evolve continuously in time.
The idea is that a well-trained ESN in the original parameter regime will maintain predictive power in the new parameter regime, provided the training of the ESN is updated with an additional amount of data from the new parameter regime. We now describe this idea mathematically.}

We describe our approach in the context of flows on metric spaces.
Let $(M,d)$ be a metric space and let $\Phi : M \times [0, \infty ) \to M$ be a continuous map that defines a flow on $M$.
That is, we have $\Phi^{t+s} = \Phi^{t} \circ \Phi^{s}$ for all $s,t \geqslant 0$, where $\Phi^{t} (\cdot ) = \Phi (\cdot , t)$.
Think of $M$ as a function space and $\Phi$ as a flow on this function space generated by an evolution partial differential equation.
Research at the interface of machine learning and dynamical systems has asked the following question.
Can $\Phi$ be learned for predictive purposes?

The single-orbit approach asks the machine learning architecture to predict individual trajectories of $\Phi$ once training is complete.
Echo state networks have been successfully used to predict trajectories of complex models, including chaotic systems~\cite{pathak2017using, estevez2023echo, hesch21, girvan2020separation, nadiga2021reservoir, chen2022predicting, doan2020physics, gonzalez2022optimizing}.
It is challenging, though, to prove rigorous mathematical lower bounds for the length of time the single-orbit predictions remain accurate.
Alternatively, one can ask the machine learning architecture to learn and then predict statistical properties of the dynamics.
This approach has also succeeded.
Indeed, echo state networks can learn statistical properties of attractors of chaotic systems~\cite{pathak2017using,pathak2018model, zhixin2018attractor,vlachas2020attractor,LI2022321,antonik2018using,racca2021robust,tanaka2019recent}.

In this paper, we focus on the following problem.
Suppose that the flow $\Phi$ depends on a parameter $\gamma$.
We write $\Phi_{\gamma}$ to indicate this dependence.
For example, $\gamma$ could be a parameter in the PDE that generates the flow $\Phi_{\gamma}$.
Suppose that for each $\gamma$ in some parameter set $\Gamma$, the flow $\Phi_{\gamma}$ admits an attractor that supports an ergodic invariant measure $\mu_{\gamma}$.
Suppose further that for each $\gamma \in \Gamma$, $\mu_{\gamma}$ describes the asymptotic distribution of the orbit of a typical $x \in M$, meaning that
\begin{equation}
\lim_{T \to \infty} \frac{1}{T} \int_{0}^{T} \delta_{\Phi_{\gamma}^{s} (x)} = \mu_{\gamma}
\end{equation}
in the weak-$\ast$ topology, \revision{where $\delta_{z}$ denotes the Dirac-$\delta$ measure at $z$.}
In this setting, we ask two questions.
First, if we fix a parameter value $\gamma_{1} \in \Gamma$ and train an echo state network using trajectories from $\Phi_{\gamma_{1}}$, can the echo state network then predict statistical properties of $(f \circ \Phi_{\gamma_{1}}^{t})_{t \geqslant 0}$, where $f : M \to \mathbb{R}^{n}$ is a given observable?
Second, assuming an affirmative answer to the first question, how might we predict the statistical properties of $(f \circ \Phi_{\gamma_{2}}^{t})_{t \geqslant 0}$ with respect to a \textit{different} parameter value $\gamma_{2} \in \Gamma$?

The point of this paper is to answer these two questions.
Our approach is as follows.
Suppose the echo state network has been trained using trajectories from $\Phi_{\gamma_{1}}$ and can successfully predict statistical properties of $(f \circ \Phi_{\gamma_{1}}^{t})_{t \geqslant 0}$.
Since invariant measures of dynamical systems can depend sensitively on parameters, the invariant measure $\mu_{\gamma_{2}}$ associated with a different parameter value $\gamma_{2}$ might differ meaningfully from $\mu_{\gamma_{1}}$.
Consequently, we would not expect the echo state network to successfully predict statistical properties of $(f \circ \Phi_{\gamma_{2}}^{t})_{t \geqslant 0}$ when trained using trajectories from $\Phi_{\gamma_{1}}$.
We conjecture that the predictive capability of the echo state network can be `transferred' from the $\gamma_{1}$-system to the $\gamma_{2}$-system using transfer learning~\cite{inubushi2020transfer,chen2022predicting}.
That is, the echo state network will successfully predict statistical properties of $(f \circ \Phi_{\gamma_{2}}^{t})_{t \geqslant 0}$ after we update its training using an \revision{additional} new dataset from the new regime. \revision{The amount of new data depends on the dynamics of the underlying equations.}

We show that our proposed approach succeeds for the generalized Kuramoto-Sivashinsky (gKS) equation, an archetypal nonlinear PDE that exhibits spatiotemporal chaos.
We train an echo state network to accurately predict the statistics of gKS dynamics.
We then show that transfer learning allows the echo state network to retain predictive power as we vary a parameter in the dispersion relation of the gKS equation or the length of the physical domain.

\revision{As discussed above, we assume that the gKS model admits an ergodic invariant measure that describes the asymptotic distribution of typical initial data.
We assume that the trained ESNs we use for prediction \textit{also} possess this property.
We have verified this hypothesis numerically.}

\revision{Although this paper focuses on long-term statistical prediction, we nevertheless show that transfer learning nontrivially improves the length of time that predictions of individual gKS trajectories remain accurate.
In particular, we report a gain of one to two Lyapunov times when using transfer learning to account for changes in the length of the spatial domain.}

Transfer learning is a popular concept in many areas of machine learning, including speech recognition and image processing.
However, transfer learning for echo state networks has received comparatively little attention in the literature.
We conclude that transfer learning can allow echo state networks to maintain predictive power \revision{when model parameters instantaneously jump from one parameter set to another.}
Such flexibility is important when statistical properties of the dynamics depend meaningfully on model parameters.


\section{KS and gKS equations}

\subsection{Background}
The one-dimensional Kuramoto-Sivashinsky (KS) equation~\cite{kuramoto1976persistent,SIVASHINSKY1977} has become a canonical model for the study of spatiotemporal chaos.
The one-dimensional KS equation is given by
\begin{equation}
\label{eqn:KS}
   \frac{\partial u}{\partial t}+\frac{\partial^2 u}{\partial x^2} + \frac{\partial^4 u}{\partial x^4}+u \frac{\partial u}{\partial x} = 0,
\end{equation}
with the periodic boundary condition $u(x+L,t) = u(x,t)$ for all $x$ and every $t \geqslant 0$.
The KS equation behaves nicely as an infinite-dimensional dynamical system.
In particular, the KS model admits a finite-dimensional attractor and is therefore in some sense equivalent to a finite-dimensional dynamical system~\cite{constantin2012integral, foias1988inertial, nicolaenko1985some}.
See~\cite{temam2012infinite} for an overview.
Moreover, the KS model admits an inertial manifold that contains the global attractor~\cite{jolly1990approximate, COLLET1993, temam1994}.
The KS equation has been studied extensively, both analytically and computationally~\cite{kevrekidis1990back, hyman1986order, smyrlis1991predicting, munkel1996intermittency}.

The parameter $L$ (the length of the spatial domain) plays the role of a bifurcation parameter for the KS model.
This parameter determines the number of linearly unstable Fourier modes, $S_{\mathrm{unst}}$. In addition, the maximum of the coefficient of the linearized equation occurs at $M_{\mathrm{unst}}$.
We refer to $M_{\mathrm{unst}}$ as ``the most unstable wavenumber".
These quantities are given by
\begin{equation}
\label{unst}
S_{\mathrm{unst}} = \left\lfloor \frac{L}{2\pi} \right \rfloor,
\qquad 
M_{\mathrm{unst}} = \frac{L}{2\sqrt{2}\pi}.
\end{equation}
Here $\lfloor \cdot \rfloor$ denotes the greatest integer function.

The one-dimensional generalized Kuramoto-Sivashinsky equation is obtained from~\eqref{eqn:KS} by adding a parametric dispersion term (third derivative).
The gKS equation is given by
\begin{equation}\label{eqn:gKS}
   \frac{\partial u}{\partial t}+\frac{\partial^2 u}{\partial x^2}+ \gamma \frac{\partial^3 u}{\partial x^3}+ \frac{\partial^4 u}{\partial x^4}+u \frac{\partial u}{\partial x} = 0,
\end{equation}
where $\gamma$ is a parameter.
The boundary conditions are again periodic.

The dynamics of the gKS equation have been studied extensively~\cite{kawahara1983formation, manneville1985liapounov, kalliadasis2011falling,  balmforth1995solitary, chang1995interaction, ei1994equation, duprat2009liquid, tseluiko2010interaction, duprat2011wave, tseluiko2014weak, pradas2011rigorous, chang1993laminarizing}.
The parameter $\gamma$ has a strong influence on gKS dynamics.
For positive values of $\gamma$ near zero, the gKS model exhibits spatiotemporal chaos, like the KS model does.
As $\gamma$ increases, the gKS model exhibits a transition toward less chaotic or even non-chaotic dynamics (see~\cite{goprka15} for a detailed numerical study).
In the limit $\gamma \to \infty$, the gKS equation is equivalent to the integrable Korteweg-de Vries (KdV) equation~\cite{tseluiko2010interaction}.

In this paper, we focus on statistical prediction as the parameters $L$ and $\gamma$ vary.


\subsection{Numerical method}

We use numerical solutions of the gKS equation to train ESNs.
In this subsection, we describe the numerical method we use to integrate the gKS equation.

We use a uniform grid on $[0,L]$ and 
apply finite differences
to discretize the derivatives in the gKS equation.
In particular, we use an explicit method to discretize nonlinear fluxes and treat linear terms implicitly.
We utilize Euler time-stepping.
Therefore, the discrete version of the gKS equation reads
\begin{equation}
\frac{u_j^{n+1} - u_j^n}{\delta t} = 
 - \frac{F_{j+1/2}^n - F_{j-1/2}^n}{\Delta x} + f(\bu^{n+1}),
\end{equation}
where $F_{j+1/2} = \frac{1}{6}((u_{j}^n)^2 + u_j^n u_{j+1}^n + (u_{j+1}^n)^2)$, $f(\bu^{n+1})$ represents the discretization of the linear second, third, and fourth derivatives, and $\bu^{n+1} = \{u_j^{n+1} : j=1,\ldots,N_x \}$ is the discrete solution vector at time $n+1$.
The second, third, and fourth derivatives are discretized using standard $O(\Delta x^2)$ central finite-difference formulas.


\section{Echo state networks and transfer learning}

Recurrent neural networks can be challenging to train.
Reservoir computing has emerged as an elegant solution to the training challenge~\cite{jaeger2001echo, jaeger2007echo, maass2002real}.
This is so because the reservoir has its own nonlinear dynamics and only a linear output layer is trained to match the observational data.
Reservoir computers can successfully predict chaotic dynamics~\cite{jaeger2004harnessing, pathak2018model, chattopadhyay2020data, bollt2021explaining, nadiga2021reservoir, vlachas2020backpropagation}.
They have been successfully applied in a variety of fields, including climate and weather prediction~\cite{nadiga2021reservoir, arcomano2020machine}, turbulent convection~\cite{pandey2020reservoir}, and others~\cite{chattopadhyay2020data, vlachas2020backpropagation}.

In this paper, we use the following standard ESN architecture.
The ESN consists of an input layer, a reservoir of $D$ neurons, and an output layer.
Connections between the neurons are encoded by an adjacency matrix $A$.
The reservoir processes successive inputs, where the reservoir state at time $t + \Delta t$, $r(t + \Delta t) \in \mathbb{R}^{D}$, is a nonlinear function of the previous reservoir state, $r(t) \in \mathbb{R}^{D}$, and the driving input, $X(t)$, given by
\begin{equation}
\label{states_1}
r(t+\Delta t) = f(A r(t) + W_{\mathrm{in}} X(t)).  
\end{equation}
Here $f = \tanh : \mathbb{R}^{D} \to \mathbb{R}^{D}$ is the element-wise activation function.
Let $N_{x}$ denote the dimension of the input.
Matrices $A \in \R^{D \times D}$ and $W_{\mathrm{in}} \in \R^{N_x \times D}$ remain constant throughout the training process.
During training, only the weights of the output layer, $W_{\mathrm{out}} \in \R^{D \times N_x}$, are tuned. 
Predictions are given by
\begin{equation}
X ( t +\Delta t) = W_{\mathrm{out}} \phi (r(t+\Delta t)),
\end{equation}
where $\phi$ is the nonlinear transformation~\cite{chattopadhyay2020data, pathak2018model} given component-wise by
\begin{equation}
  \tilde{r}_j (t) = \phi (r_j (t))=
    \begin{cases}
      r_j^2(t) & \text{if $j$ is even}\\
      r_j (t) & \text{if $j$ is odd}.
    \end{cases}       
\end{equation}
We found in our experiments that this nonlinear transformation slightly improves the quality of predictions for individual trajectories.

The input matrix $W_{\mathrm{in}}$ is generated component-wise from a uniform distribution on $[-\beta_{1}, \beta_{1}]$.
The adjacency matrix is generated as $A = \beta_2 W_0 / \lambda$, where $W_0$ is a sparse matrix (typically turning on less than $10\%$ of the possible connections), $\lambda$ denotes the largest eigenvalue of $W_0$, and $\beta_2$ is a scaling parameter.

The weights of $W_{\mathrm{out}} \in \R^{{D} \times N_x}$ are determined through a linear regression with $L^2$-regularization by
\begin{equation}
W_{\mathrm{out}} = \argmin_{W} \|W\mathbf{R}-\mathbf{X}\|_2^2 + \mu \|W \|_2^2,
\label{W_out_optimization_formula}
\end{equation}
where \revision{$\mu$ is a regularization parameter}, 
$\mathbf{X}$ is the training data (trajectories), and $\mathbf{R}$ denotes the reservoir timeseries (with $A$ and $W_{\mathrm{in}}$ fixed).
$W_{\mathrm{out}}$ can be expressed explicitly as
\begin{equation}
W_{\mathrm{out}} = \X \rr' ( \rr\rr' +\mu I)^{-1}.
\label{W_out_formula}
\end{equation}
Here, $\rr'$ denotes the transpose of $\rr$. 
Throughout this paper, we use ESN parameters $D=5000$
(same as in \cite{pathak2018model}), $\beta_1 = 0.01$, $\beta_2 = 0.1$, and $\mu =0.000005$, unless stated otherwise.

\revision{The ESN setup we use in this paper is similar to that of~\cite{chen2022predicting}.
The common approach is that the values of $\beta_1$ and $\beta_2$ should be 
\revisiontwo{chosen so that $A r(t) + W_{\mathrm{in}} X(t)$ evolves in the nonlinear regime of $\tanh (\cdot )$} and the ESN operates in the \revisionthree{unsaturated regime.}
\revisiontwo{Choosing appropriately depends upon the problem at hand.}
We performed a basic search to determine optimal hyperparameter values.
Our search revealed that the ESN performs considerably suboptimally when $\beta_2 \lesssim 1$ and $\beta_1 \in [0.1,0.5]$.
The performance of the ESN is not very sensitive to hyperparameter values as long as $\beta_1 \in [0.01, 0.1]$ and $\beta_2 \in [0.1, 0.4]$.
The value of $\mu$ has been taken from~\cite{chattopadhyay2020data}.}
\revisiontwo{Since our study involves long stationary simulations where initial transient behavior is not important, we do not use reservoir warm-up in this work.
We performed a synchronization test (not included here for the sake of brevity), showing that the reservoir synchronizes with 2 to 4 timesteps.}

Transfer learning~\cite{inubushi2020transfer, chen2022predicting} is a machine-learning approach wherein knowledge in a source domain is used to improve predictive performance in a target domain.
In our context, the source and target domains are simply two sets of gKS parameters.
Our central premise for linking transfer learning to ESNs is the following.
Once the ESN is trained well in the source domain, it will retain predictive power in the target domain after receiving a training update using a \revision{relatively} small transfer-learning dataset from the target domain. \revision{However, the amount of transfer learning data might depend on a particular dynamic regime.}

Writing $\X_{\mathrm{TL}}$ for the transfer-learning dataset, the update to the output layer is obtained by solving the optimization problem given by
\begin{equation}
\delta W = \argmin_{\delta W} \left\| (W_{\mathrm{out}}+\delta W) \rr_{\mathrm{TL}} - \X_{\mathrm{TL}} \right\|_2^2 + \alpha \| \delta W \|_2^2 .
\end{equation}
Here, $\alpha$ is the rate of transfer learning and $\rr_{\mathrm{TL}}$ is the (time-dependent) state of the reservoir (with the same number of time-snapshots as $\X_{\mathrm{TL}}$).
The update to the output layer and new output layer are given explicitly by
\begin{equation}
\label{transfer_learning}
\begin{aligned}
    \delta W  &= \left(\X_{\mathrm{TL}} \rr_{\mathrm{TL}}' - W_{\mathrm{out}} \rr_{\mathrm{TL}} \rr'_{\mathrm{TL}}  \right)
    \left(\rr_{\mathrm{TL}} \rr'_{\mathrm{TL}}+ \alpha I \right)^{-1} , \\
    W_{\mathrm{out}}^{\mathrm{new}} &= W_{\mathrm{out}} + \delta W.
\end{aligned}
\end{equation}
Here, $W_{\mathrm{out}}$ and $W_{\mathrm{out}}^{\mathrm{new}}$ are the output layer in the source and target domains, respectively.

\revision{In this paper, we use $\alpha=0.005$.
We performed a grid search to determine the optimal value of this hyperparameter.
In particular, we tested several values of $\alpha$, including $0.1$, $0.01$, $0.005$, $0.002$, and $0.001$.
Empirically, $\alpha = 0.005$ works well.}

Transfer learning has been used with ESNs to successfully predict individual trajectories for shallow-water equations~\cite{chen2022predicting, shahalam2024} and the gKS equation~\cite{shahalam2024}. 
In this paper, we show that transfer learning allows us to successfully track changes in the statistical signature of the gKS equation \revision{when either the size of the spatial domain or the dispersion relation instantaneously jumps from one regime to another.}


\section{Results}
\label{sec:results}

We focus on an important statistical feature of gKS dynamics, namely the power spectrum.
The \revision{\textit{time-averaged}} power spectrum is defined as the map from $\mathbb{Z}$ into $\mathbb{R}$ given by
\begin{equation}
\label{ek}
k \mapsto e_k = \frac{1}{T} \int\limits_0^T |\uh_k(t)|^2 \, \mathrm{d}t,
\end{equation}
where $\hat{u}_{k} (t)$ is coefficient $k$ at time $t$ of the Fourier-series representation of the solution to the gKS equation,
\begin{equation}
u(x,t) = \sum_{k \in \mathbb{Z}} \hat{u}_{k} (t) e^{(2 \pi i k / L) x}.
\end{equation}
The quantity $e_{k}$ represents the average energy in Fourier wavenumber $k$.
Since the solution $u(x,t)$ of the gKS equation is real-valued, the Fourier coefficients satisfy $\uh_{-k} = \uh_k^*$.
We therefore consider only nonnegative wavenumbers.

Throughout Section~\ref{sec:results}, we work with the gKS equation in the chaotic regime ($\gamma = 0$ or $\gamma$ is small and positive).
We assume the gKS equation admits an attractor that supports an ergodic measure.
We assume that this ergodic measure describes the distribution of the orbit of the initial data we select.
(As $\gamma$ increases, the gKS equation eventually becomes less chaotic~\cite{goprka15}.
In particular, the power spectrum depends on the initial data when $\gamma$ is sufficiently large.)

We first consider the KS equation and we vary the parameter $L$.
Varying $L$ is equivalent to varying the number of linearly unstable Fourier modes.
\revision{For each of several values of $L$, we train the ESN using trajectories of the KS equation with spatial domain $[0,L]$.
We show that the trained ESN accurately predicts the power spectrum of the KS equation when the spatial domains used for training and prediction are the same.}

We then show that transfer learning is effective for various tuples $(L_{1}, L_{2})$ of length parameters.
That is, the ESN accurately predicts the power spectrum of the KS equation with parameter $L_{2}$ \revision{and the correlation dimension of the corresponding attractor} after a two-step training process.
First, the ESN is trained using trajectories of the KS equation with parameter $L_{1}$.
Second, the training is updated using \revision{some} trajectory information from the KS equation with parameter $L_{2}$.

\revision{We complete our suite of experiments for the KS equation by showing that transfer learning (with respect to $L$) nontrivially improves the length of time that predictions of individual KS trajectories remain accurate.
In particular, we report a gain of one to two Lyapunov times using a modest percentage of data for transfer learning ($5\%$).}

We conclude Section~\ref{sec:results} by carrying out the program described above for the gKS equation.
This time, we fix $L$ and vary the parameter $\gamma$ in the dispersion relation.
\revision{We show that when $\gamma$ is small, the well-trained ESN accurately predicts the power spectrum and the correlation dimension of the attractor, when the values of $\gamma$ used for training and prediction match.
Further, \revisiontwo{we show that} transfer learning is effective in the small-$\gamma$ regime.}

\textbf{Parameter selection for experiments involving the power spectrum.}
\revision{In the previous study~\cite[Section~5.1.3]{shahalam2024}, we determined that $N_x=256$, averaging time $T=10000$ (approximately $1000$ Lyapunov times), integration timestep $\delta t=0.001$, and snapshot sampling timestep $\Delta t=0.25$ are sufficient for resolving power spectra.
We use these parameters for the experiments in this paper, unless specifically stated otherwise.
Since we always use $N_x=256$ points in for the discretization, $\Delta x$ depends on $L$.}


\subsection{Predicting the statistical behavior of the KS equation without transfer learning}
\label{SubSec:KSnoTL}

We show that for each of several values of $L$, the ESN accurately predicts the power spectrum of the KS equation with spatial domain $[0,L]$, once trained using KS trajectories from the same spatial domain.
We consider four values of $L$ that correspond to four different numbers of linearly unstable Fourier modes.
Table~\ref{tab:L_data} lists each choice of $L$, together with the corresponding number of linearly unstable Fourier modes ($S_{\mathrm{unst}}$), the most unstable wavenumber ($M_{\mathrm{unst}}$), and an approximation of $\lambda_{\max}$, the maximal Lyapunov exponent.
We use the approximation for the Lyapunov exponents ($\lambda_{0} \geqslant \lambda_{1} \geqslant \cdots \geqslant \lambda_{i} \geqslant \cdots$) given in~\cite{edson2019lyapunov} by
\begin{equation}
\label{Lyapunov-exponent-formula}
\lambda_{i} (L) \approx 0.093 - \frac{(0.94) (i-0.39)}{L}.
\end{equation}

\begin{table}[ht]
\centering
\caption{
Values of $L$ for which we predict the power spectrum of the KS equation with spatial domain $[0,L]$.
For each value of $L$, we list the number of linearly unstable Fourier modes ($S_{\mathrm{unst}}$) given by~\eqref{unst}, the most unstable wavenumber ($M_{\mathrm{unst}}$) given by~\eqref{unst}, and an approximation of the maximal Lyapunov exponent (given in~\cite{edson2019lyapunov} and by~\eqref{Lyapunov-exponent-formula}).
}
\label{tab:L_data}
\begin{tabular}{| c | c | c | c |}
\hline
L & $S_{\mathrm{unst}}$ & $M_{\mathrm{unst}}$ & $\lambda_{\max}$\\
\hline
22 & 3 & 2.48 & 0.1097 \\
\hline
29 & 4 & 3.26 & 0.1056 \\
\hline
35 & 5 & 3.94 & 0.1035\\
\hline
43 & 6 & 4.84 & 0.1015 \\
\hline
\end{tabular}
\end{table}

We performed the following experiment for each value of $L$.
We began with an untrained ESN.
We generated $30$ KS trajectories of length \revision{$T = 10000$ (approximately $1000$ Lyapunov times)} using random initial data of the form
\begin{equation}\label{KSIC}
    u(x,0)=c_1 \cos\left(\frac{p_1x\pi}{L}\right)+c_2 \cos\left(\frac{p_2x\pi}{L}\right) ,
\end{equation}
where we independently selected $c_1 , c_2 \sim \mathrm{Uniform} ([0,1])$ and $p_1 , p_2 \sim \mathrm{Uniform} ( \{ 1,2,3,4,5,6 \} )$.
We used $20$ of these trajectories to train the ESN and the other $10$ for prediction.

\begin{figure}[ht]
\centering
\centerline{
    \includegraphics[width=0.48\textwidth]{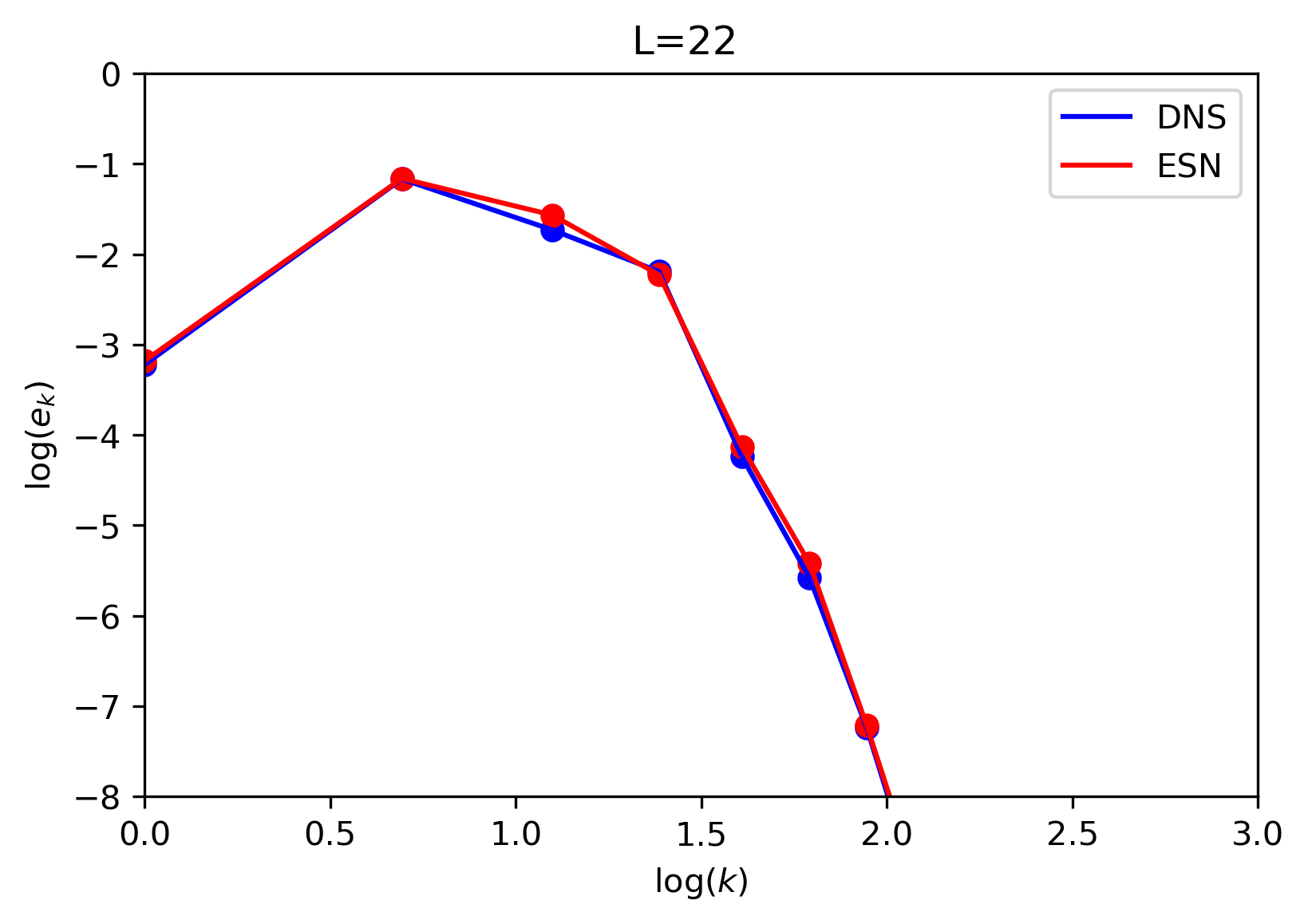}
    \includegraphics[width=0.48\textwidth]{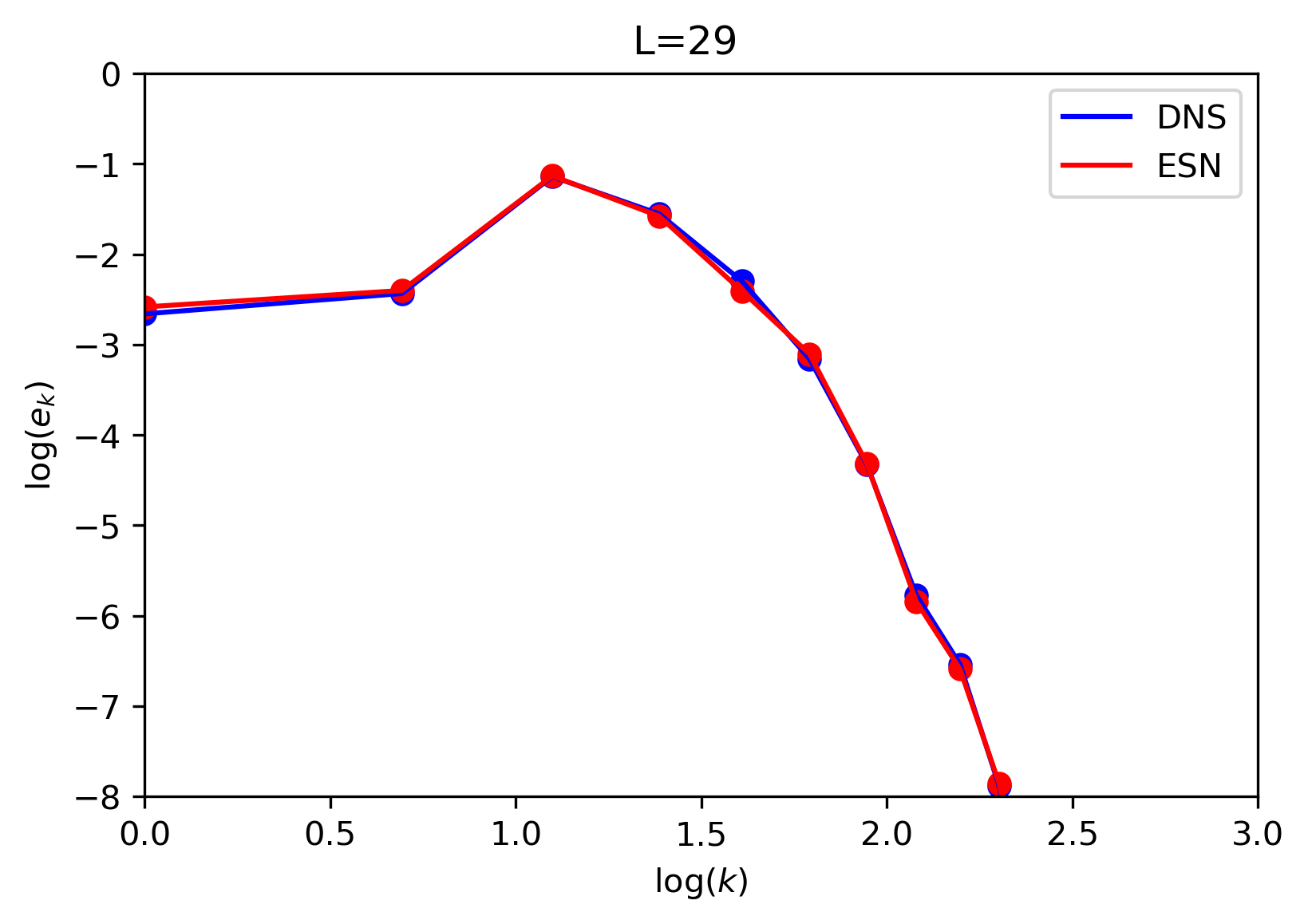}
  }
  \centerline{
    \includegraphics[width=0.48\textwidth]{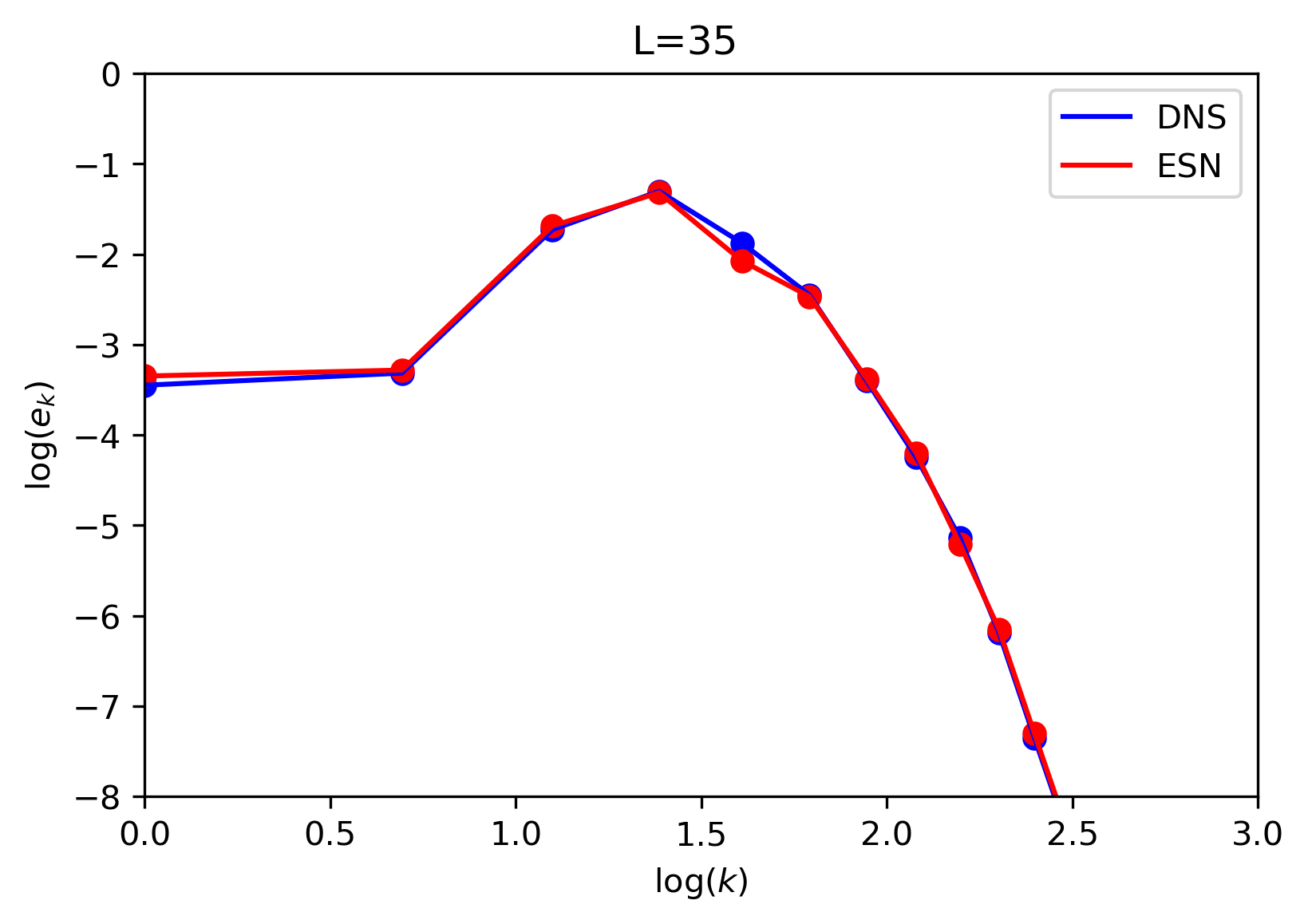}
    \includegraphics[width=0.48\textwidth]{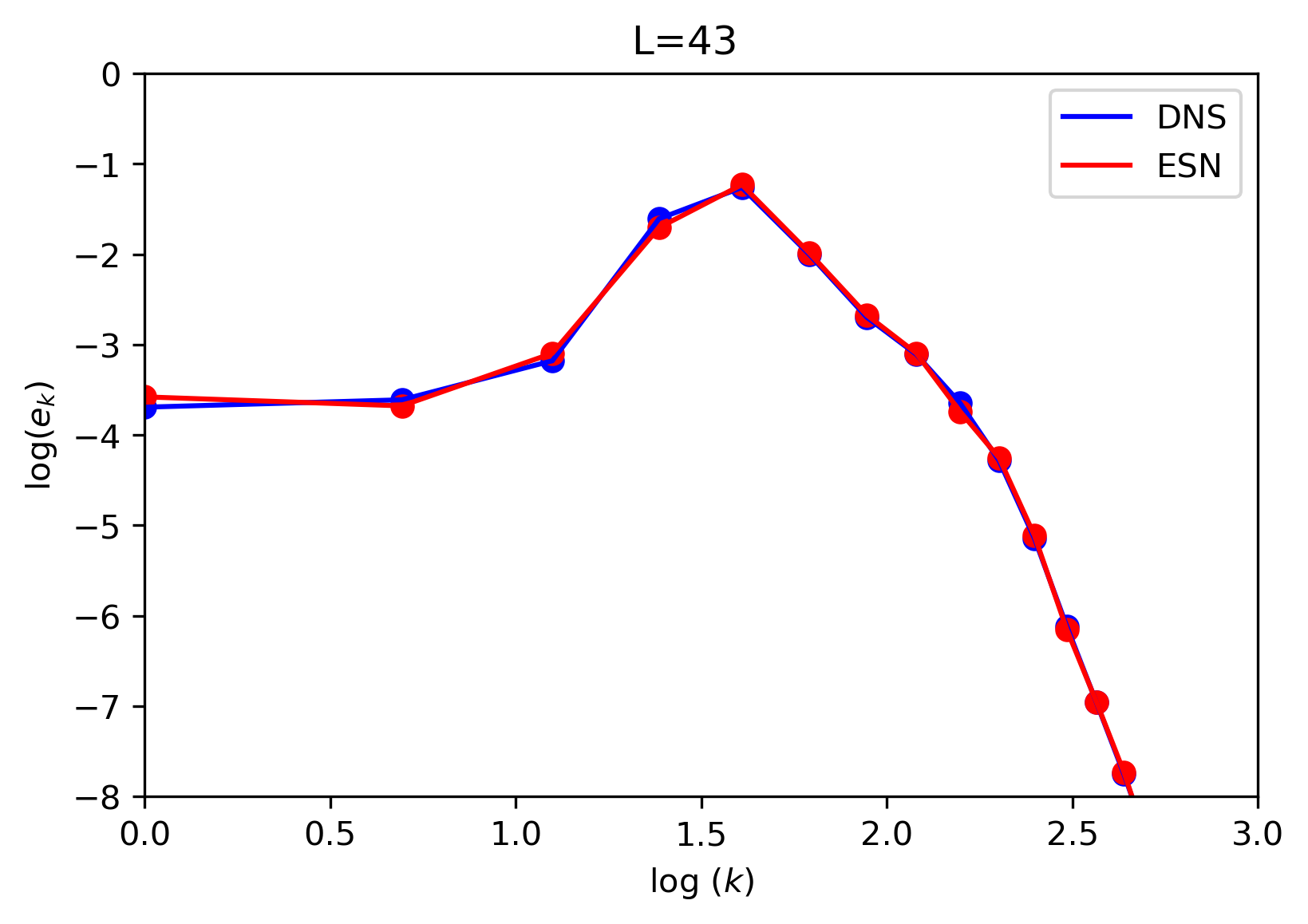}
}
\caption{Once trained on KS trajectories from the spatial domain $[0,L]$, the ESN then accurately predicts the \revision{\textit{time-averaged}} power spectrum of the KS equation with the same spatial domain.
Each blue curve shows the power spectrum obtained from direct numerical simulation, averaged over $10$ DNS trajectories used for prediction.
Red curves show the ESN predictions of the power spectrum.
Importantly, the ESN predictions of the power spectrum remain accurate long after ESN predictions of individual trajectories become inaccurate.
Spatial domain sizes $L=22$ (top left), $L=29$ (top right), $L=35$ (bottom left), and $L=43$ (bottom right) correspond to four different numbers of linearly unstable Fourier modes (see Table~\ref{tab:L_data}).
Plots show $\log (e_{k})$ versus $\log (k)$, where $k$ denotes wavenumber.
}
\label{fig:L_no_TL}
\end{figure}

Figure~\ref{fig:L_no_TL} demonstrates that for each of the four values of $L$ in Table~\ref{tab:L_data}, the ESN accurately predicts the power spectrum of the KS equation (once trained).
Each panel compares the power spectrum obtained by direct numerical simulation (DNS) to the ESN prediction.
We average over $10$ DNS trajectories used for prediction to obtain each DNS-based curve.

\subsection{Predicting the statistical behavior of the KS equation when the length of the spatial domain instantaneously increases using transfer learning}
\label{sec:KSTL}

When $L$ \revision{changes instantaneously}, it is computationally expensive to train a \textit{de novo} ESN to predict the power spectrum for each new value of $L$ \revision{because one needs to generate 
a new training dataset for each new such value.
Thus, in the context of using ESNs for systems with different values of parameters, generating new training data can be a computational bottleneck.
Transfer learning can potentially reduce the computational cost by utilizing smaller datasets for retraining.}
\revision{In this section, we demonstrate} the efficacy of transfer learning with respect to $L$ for the prediction of the power spectrum.
\revision{Further, we show that when $L$ increases instantaneously, transfer learning restores the ability of the ESN to predict the correlation dimension of the attractor.}

\begin{figure}[ht]
\centering
\includegraphics[width=0.32\textwidth]{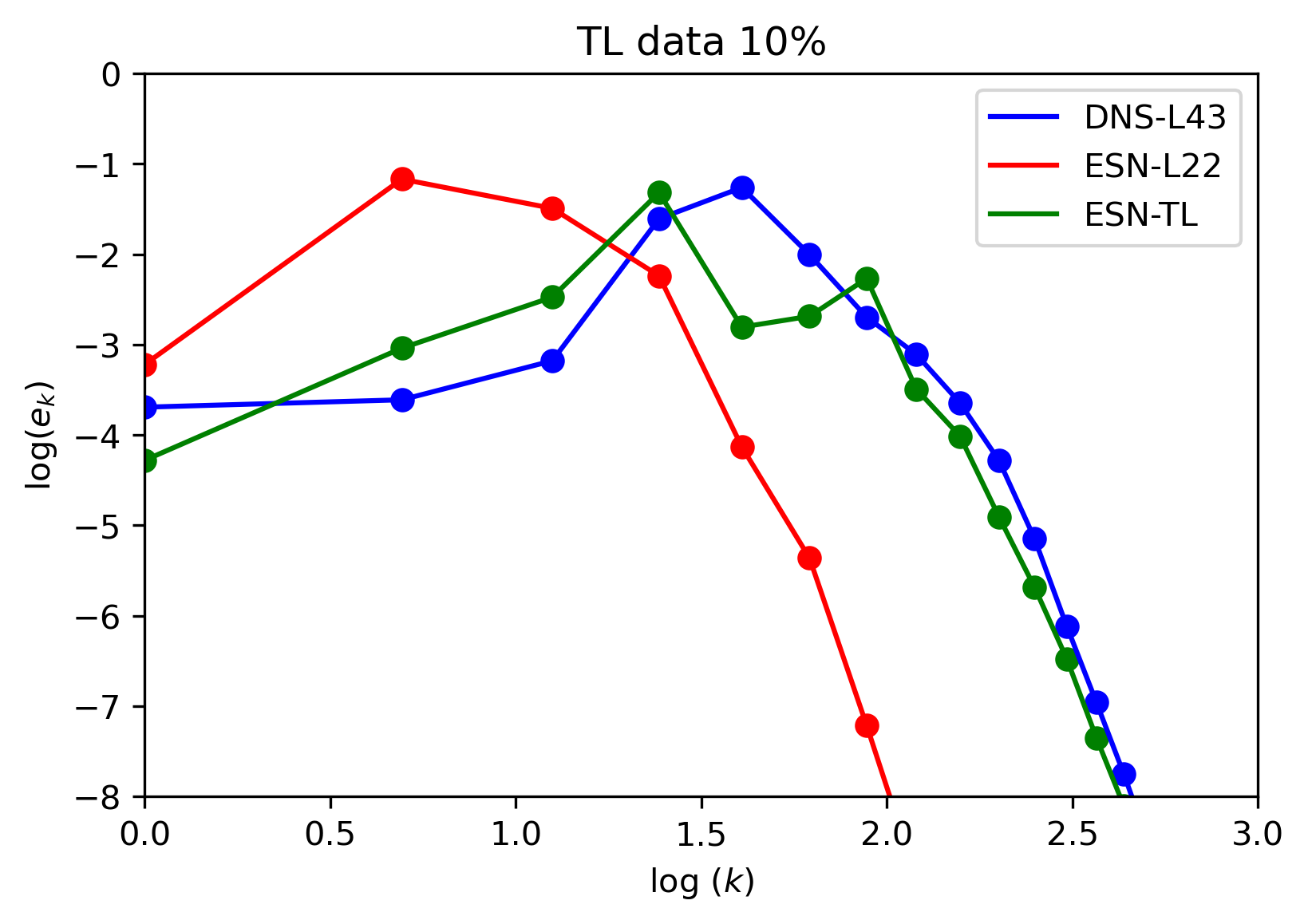}
\includegraphics[width=0.32\textwidth]{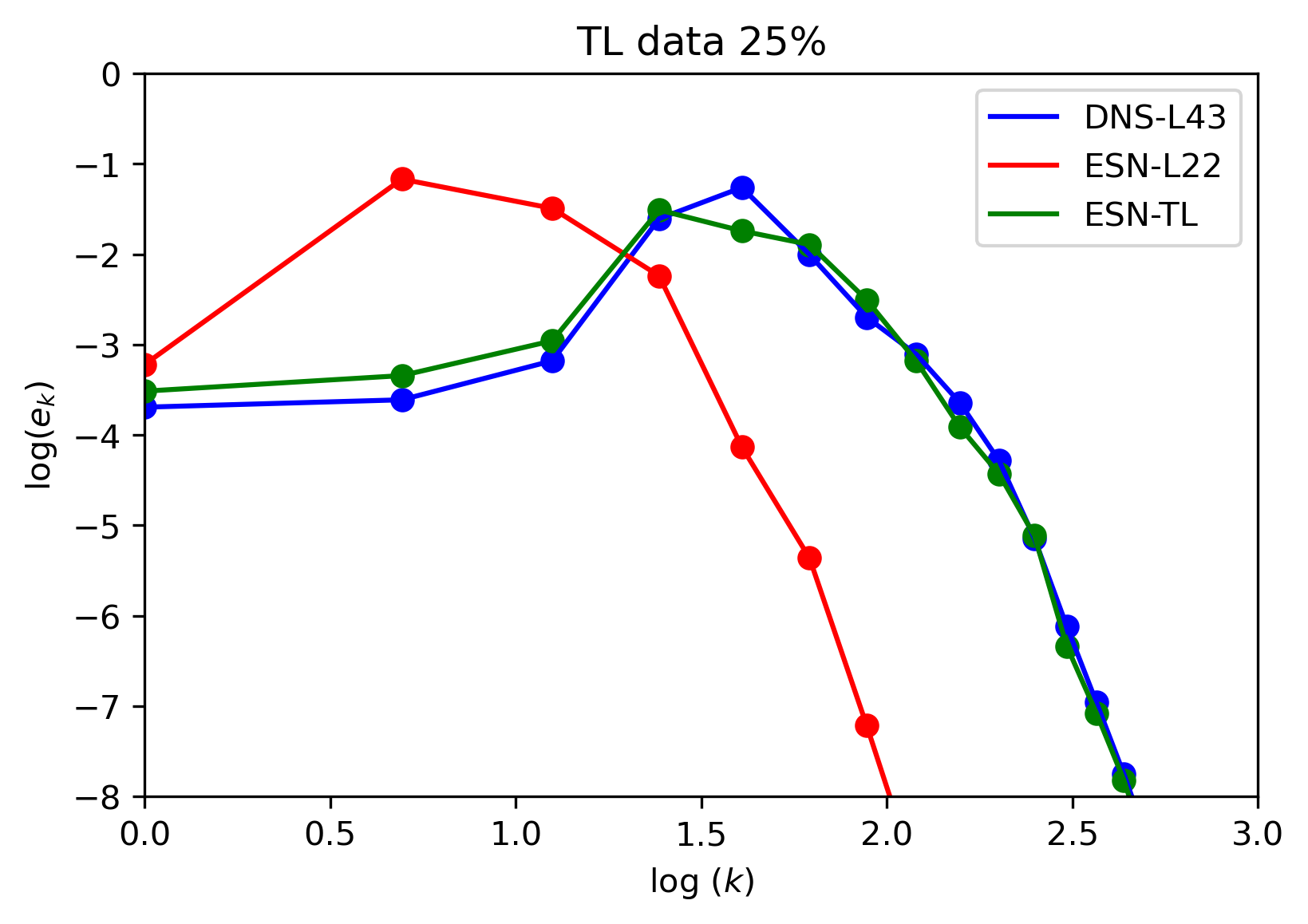}
\includegraphics[width=0.32\textwidth]{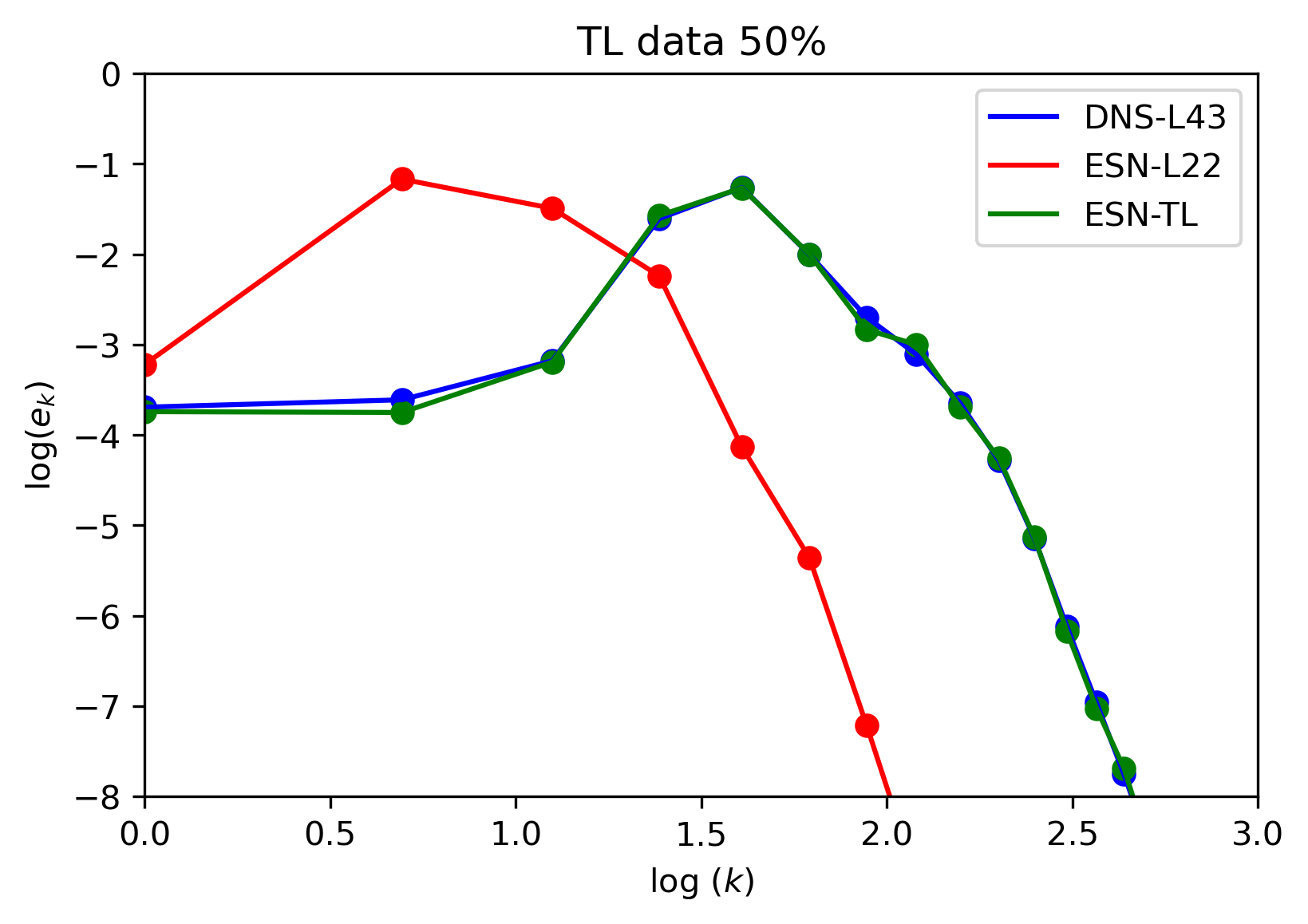}
\caption{
Transfer learning allows the ESN to accurately predict the \revision{\textit{time-averaged}}
power spectrum of the KS equation as the size of the spatial domain varies.
Curves illustrate power spectra.
DNS-L43: Power spectrum obtained by averaging over the $10$ DNS trajectories in the $L=43$ regime used for prediction.
ESN-L22: Power spectrum for the $L=22$ regime predicted by the ESN, once it is trained on $20$ trajectories from the $L=22$ regime.
ESN-TL: Power spectrum for the $L=43$ regime predicted by the ESN after transfer learning of level $10\%$ (left), $25\%$ (middle), and $50\%$ (right).
Plots show $\log (e_{k})$ versus $\log (k)$, where $k$ denotes wavenumber.
}
\label{fig:L_with_TL}
\end{figure}

We consider $L=22$ and $L=43$, values of $L$ that correspond to the presence of $3$ and $6$ linearly unstable Fourier modes, respectively.
Importantly, the power spectrum for $L=43$ differs considerably from the power spectrum for $L=22$ (Figure~\ref{fig:L_with_TL}, blue and red curves, respectively).

We performed the following experiment.
First, we trained the ESN in the $L=22$ regime using $20$ trajectories from that regime.
Second, we set the transfer learning rate to $\alpha = 0.005$ and then updated the training of the ESN using trajectory data from the $L=43$ regime.
We report the amount of TL data used for the training update \revision{relative to the size of the dataset used to train the ESN when $L=43$ in the experiment from Section~\ref{SubSec:KSnoTL} ($20$ trajectories).}

In the absence of transfer learning, the ESN (trained only on trajectories from the $L=22$ regime) \revision{poorly predicts the power spectrum in the $L=43$ regime.}
However, transfer learning improves the prediction of the power spectrum in the $L=43$ regime. 
Figure~\ref{fig:L_with_TL} illustrates this improvement as a function of amount of TL data.
The blue curve (same curve in each panel) is the power spectrum in the $L=43$ regime, computed by averaging over the $10$ DNS trajectories from this regime used for prediction.
This curve serves as ground truth.
The red curve (same curve in each panel) is the power spectrum for $L=22$ predicted by the ESN after having been trained on $20$ DNS trajectories from the $L=22$ regime.
Comparing the red and blue curves, we see that the power spectrum for $L=22$ differs considerably from that for $L=43$.
The green curves illustrate the ESN prediction of the power spectrum for the $L=43$ regime after transfer learning of level $10\%$ (left), $25\%$ (center), and $50\%$ (right).
The quality of the prediction based on transfer learning rapidly improves as the level of TL data increases.
In particular, the prediction is \revisiontwo{excellent} at level $50\%$ (right).

\begin{figure}[ht]
\centering
\centerline{
\includegraphics[width=0.7\linewidth]{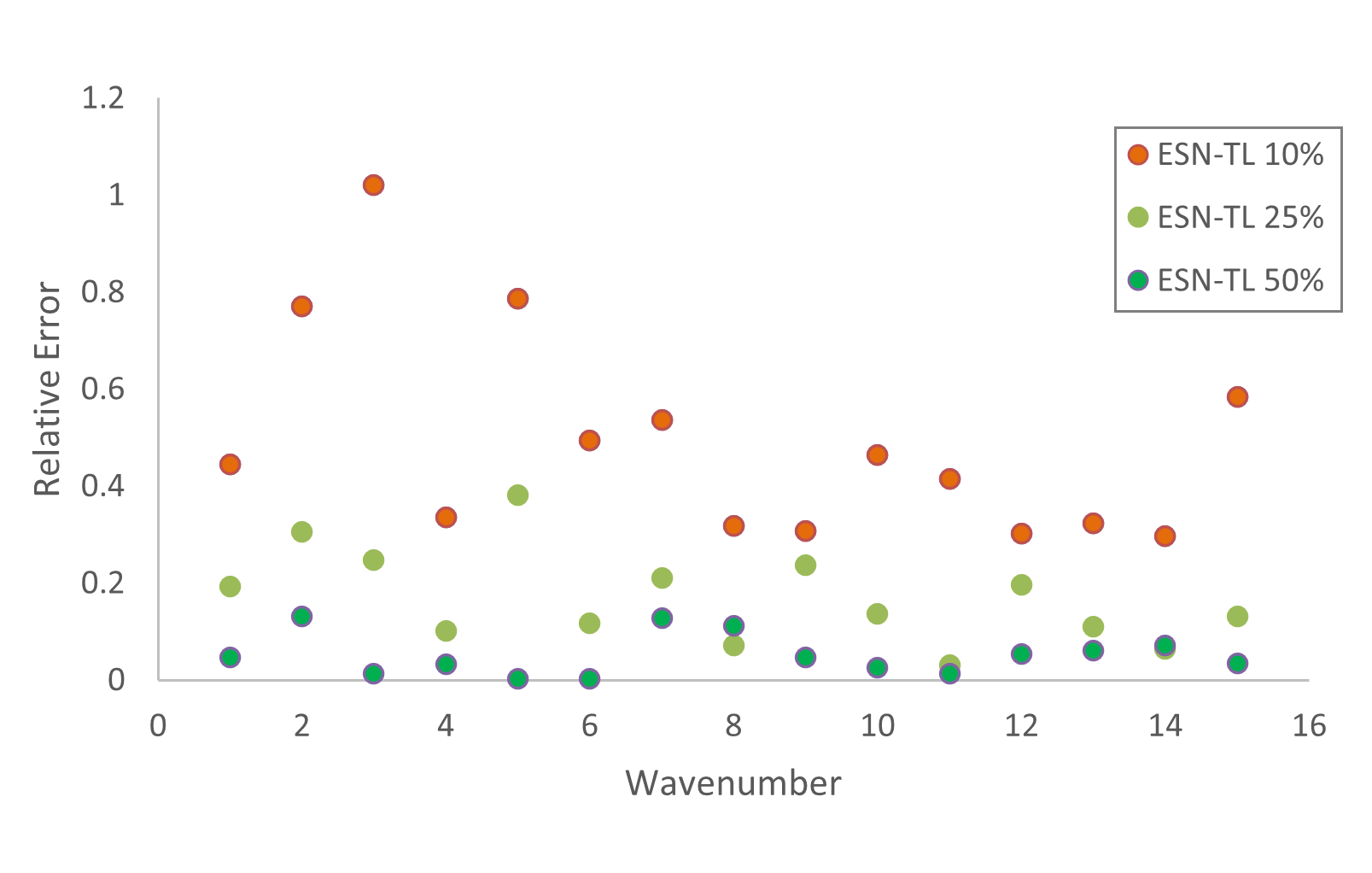}}
\caption{Plot corresponding to Table~\ref{tab:error_L_TL}.
}
\label{fig:scatterplot_L}
\end{figure}

\revision{We quantify the quality of the TL-based predictions of the power spectrum for $L=43$ by reporting the relative error
\begin{equation}
\label{e:relative_error}
\frac{|e_{k}^{\mathrm{DNS}} - e_{k}^{\mathrm{ESN}}|}{e_{k}^{\mathrm{DNS}}}
\end{equation}
associated with each of the first $15$ wavenumbers.
See Table~\ref{tab:error_L_TL} \revisiontwo{and Figure~\ref{fig:scatterplot_L}} for this error data.
Note the rapid improvement in the relative error associated with wavenumber $5$, the wavenumber that contains a plurality of the energy.}

\revision{Transfer learning restores the ability of the ESN to predict structural features of the attractor as well.
We consider one such feature in the paper, the correlation dimension of the attractor.
Table~\ref{tab:L_TL_corrDim} reports correlation dimension, computed using~\cite{scikitdim21}.
From left to right, the columns correspond to the red curve, then the three green curves, and finally the blue curve in Figure~\ref{fig:L_with_TL}.
Observe that the TL-based predictions of the correlation dimension approach ground truth (correlation dimension for $L=43$ computed using DNS data) as the amount of TL data increases.}

\begin{table}[ht]
\centering
\caption{
\revision{Transfer learning restores the ability of the ESN to predict the correlation dimension of the attractor when the length of the spatial domain instantaneously increases.
From left to right, the columns correspond to the red curve, then the three green curves, and finally the blue curve in Figure~\ref{fig:L_with_TL}.}
}
\label{tab:L_TL_corrDim}
\begin{tabular}{| c | c | c | c | c |}
\hline
ESN $L=22$ &
ESN TL $10\%$ & ESN TL $25\%$ & ESN TL $50\%$ & DNS $L=43$ \\
\hline
3.721 & 5.418 & 5.736 & 6.357 & 6.442 \\
\hline
\end{tabular}
\end{table}

\revision{\textbf{An ergodicity check.}
We performed a variant of the experiment described in Figure~\ref{fig:L_with_TL} in order to empirically verify ergodicity.
This time, we changed the $20$ trajectories we used for training in the $L=22$ regime and we used one long DNS trajectory ($T=40000$) for prediction.
Table~\ref{tab:error_ergodicity} reports relative error as a function of wavenumber for this ergodicity experiment.
As expected, relative errors in Tables~\ref{tab:error_L_TL} and~\ref{tab:error_ergodicity} closely match for those wavenumbers that contain significant energy.}


\subsection{Transfer learning upgrades single-orbit prediction for the KS model}
\label{sec:single_orbit}

\revision{Transfer learning has been used with ESNs to successfully predict individual trajectories for shallow-water equations~\cite{chen2022predicting, shahalam2024}.
In this section, we show that transfer learning meaningfully increases the length of time that single-orbit predictions for the KS model remain valid.}

\revision{Suppose an ESN is trained on KS trajectories that correspond to spatial domain length $L=22$.
How well would this ESN then predict single-orbit dynamics for larger values of $L$?
Alam~\cite[Section~4.3.2]{shahalam2024} found that the timescale for which predictions of individual trajectories remain accurate decays approximately linearly with $L$ in this scenario.
If $L > 22$, one can potentially deal with this decay by training a fresh reservoir computing architecture on $L > 22$ trajectories.
Architecture options here include a larger ESN or a coupled network of ESNs.
Transfer learning provides another potential solution.
In this section, we show that transfer learning meaningfully enhances single-orbit predictive performance for the KS model when $L > 22$.}

\revision{We performed the following experiment.
We initially trained the ESN using $20$ trajectories (each of length $T=20000$, or approximately $2000$ Lyapunov times) using $20$ KS trajectories from the $L=22$ regime.
We used this particular ESN (no transfer learning) to then predict $5$ individual KS trajectories for integer increments of $L$ from $L=23$ to $L=29$.
Table~\ref{tab:single_orbit} (top) reports times (measured in Lyapunov times) for which relative prediction errors in $L^{2}$ remain below $10\%$.
As expected, we observe decay as $L$ increases.}

\revision{To assess the efficacy of transfer learning, we performed the following experiment for each integer value of $L$ from $L=23$ to $L=29$.
We updated the training of the original ESN using one trajectory of length $20000$ from the new $L$ regime (a TL level of $5\%$ relative to the size of the original training dataset) and then used the new ESN to predict the same $5$ KS trajectories.
Table~\ref{tab:single_orbit} (bottom) reports times (measured in Lyapunov times) for which relative prediction errors in $L^{2}$ remain below $10\%$ for the TL-based predictions.
Predictive performance typically improves by $1$ to $2$ Lyapunov times (bold indicates strict increase).}

\begin{figure}[ht]
\centering
\centerline{
    \includegraphics[width=0.48\textwidth]{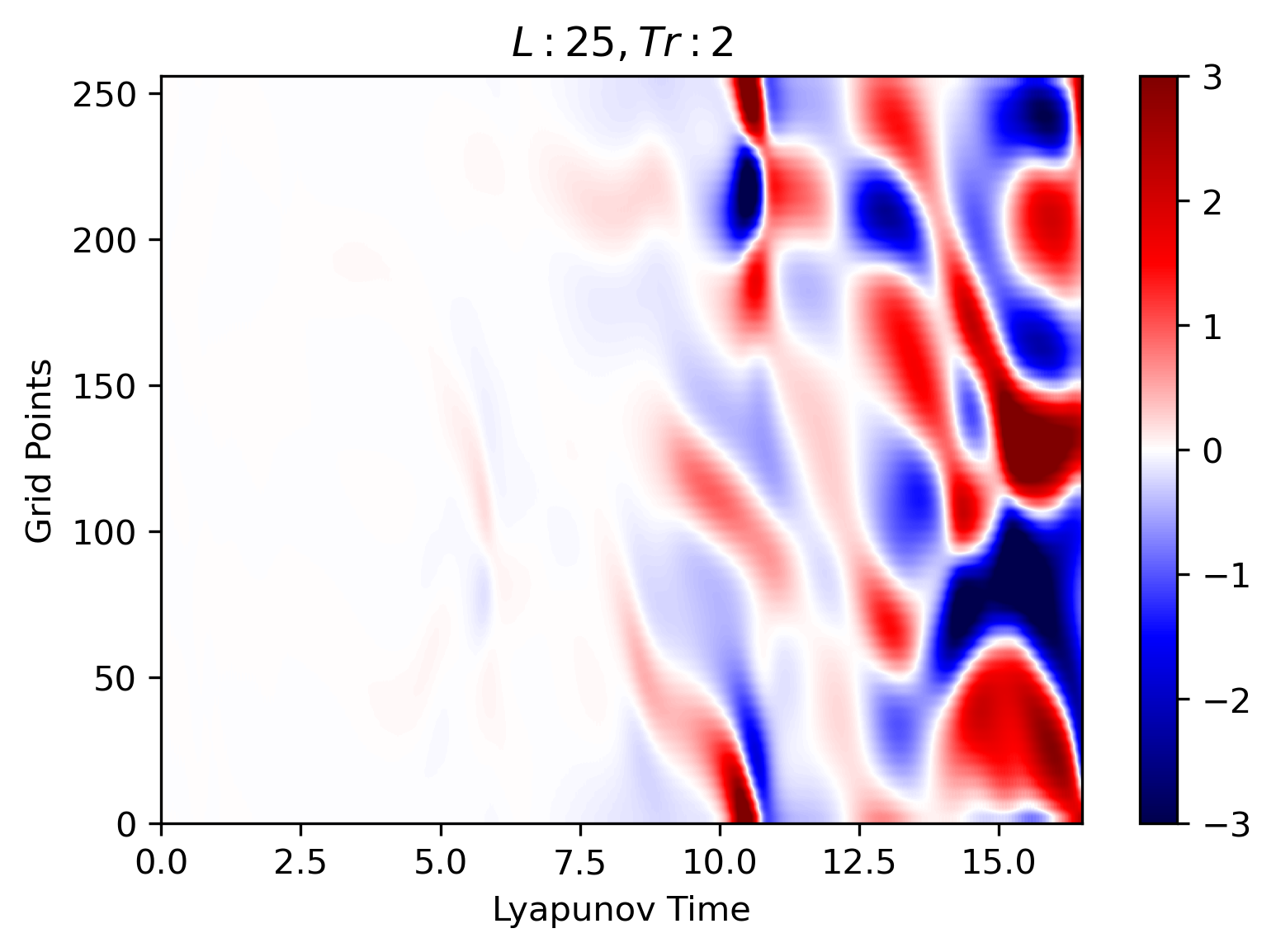}
    \includegraphics[width=0.48\textwidth]{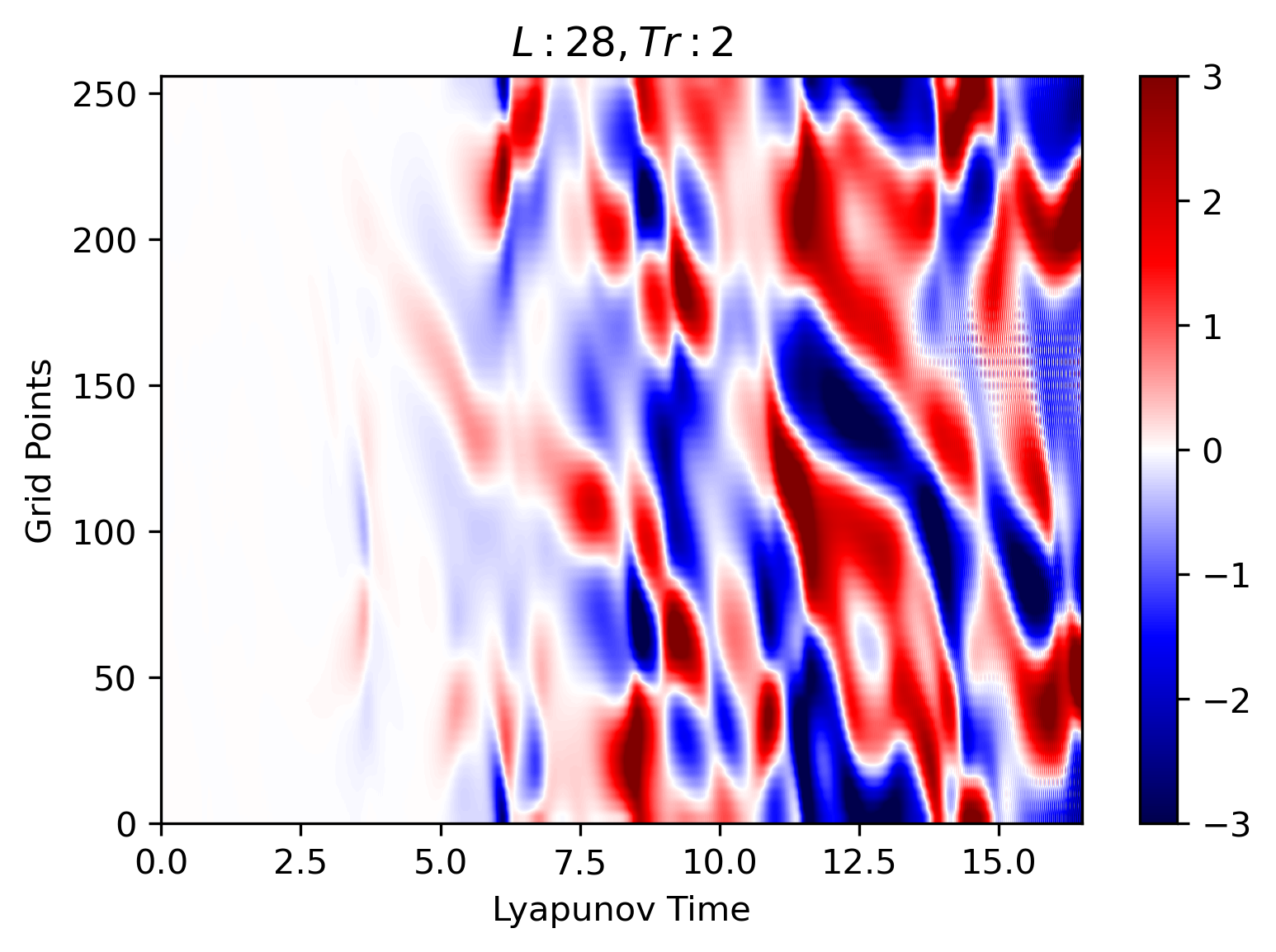}
}
  \centerline{
    \includegraphics[width=0.48\textwidth]{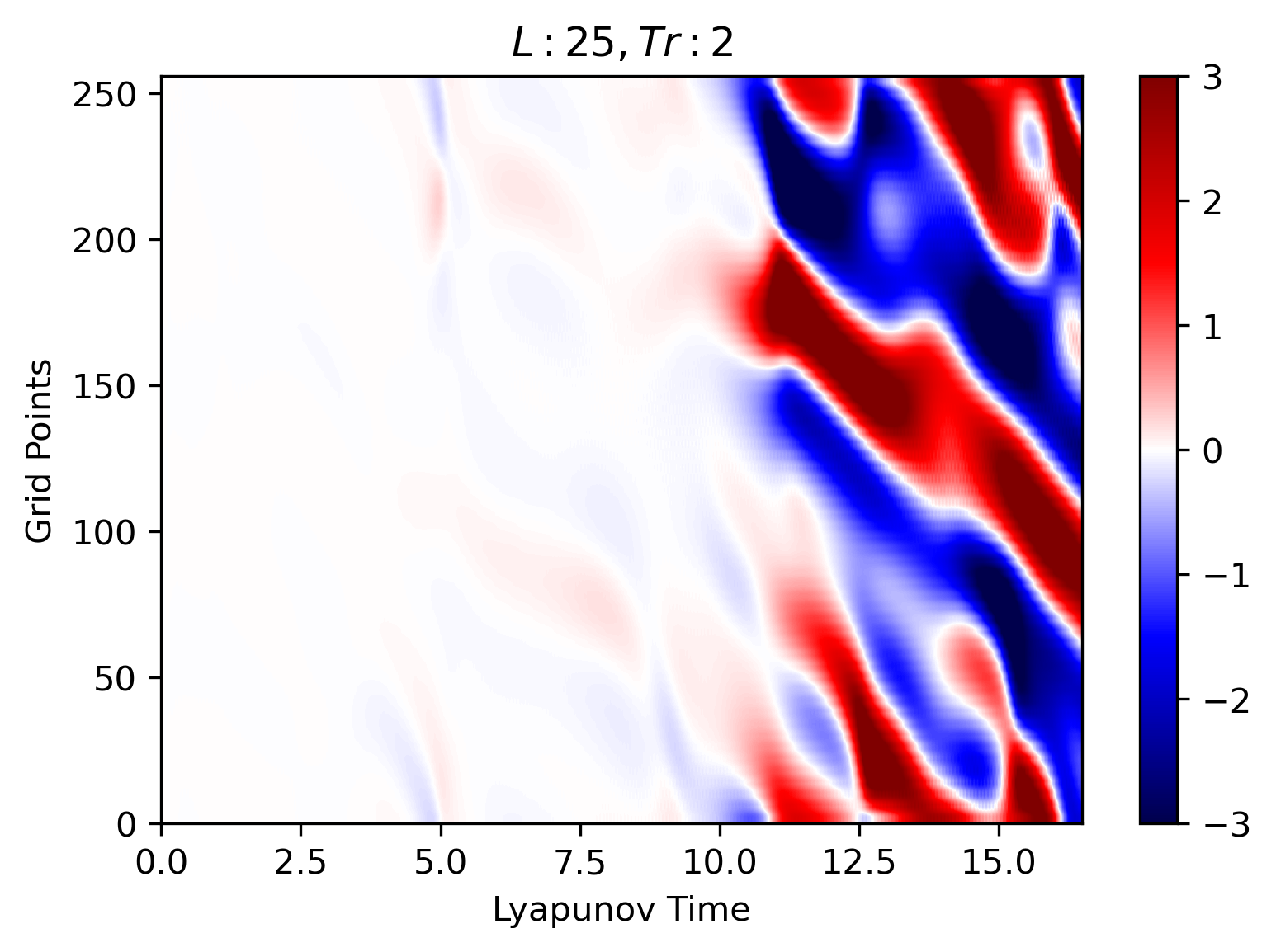}
    \includegraphics[width=0.48\textwidth]{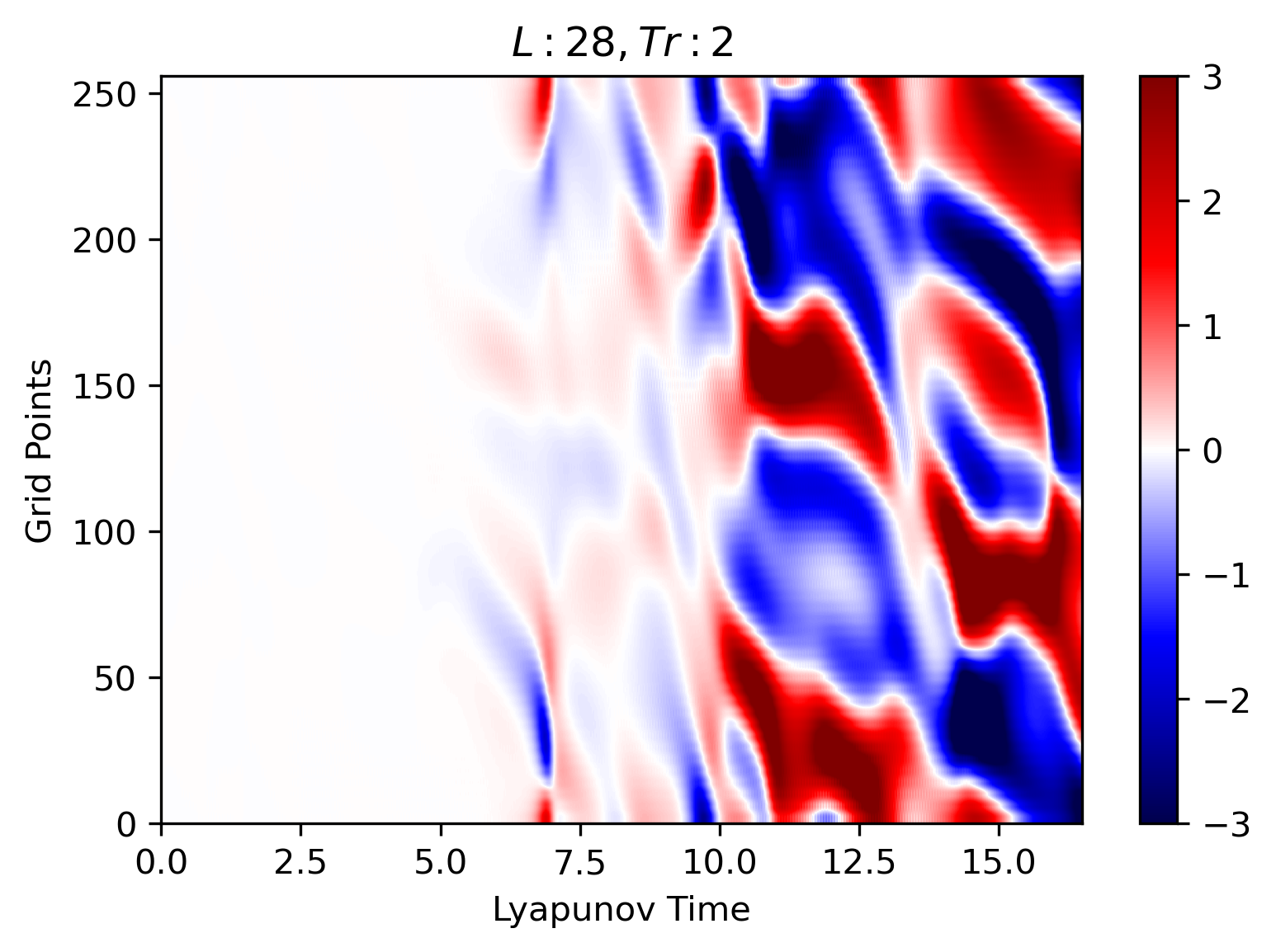}
}
\caption{\revision{Difference between the ESN-predicted KS trajectory and the KS trajectory computed by direct numerical simulation without transfer learning (top row) and with transfer learning (bottom row).
Left and right columns correspond to Trajectory~2 for $L=25$ and $L=28$, respectively, in Table~\ref{tab:single_orbit}.}
\label{fig:single_orbit}
}
\end{figure}

\revision{Figure~\ref{fig:single_orbit} illustrates the difference between the ESN-predicted KS trajectory and the DNS-computed KS trajectory ($\mathrm{ESN}-\mathrm{DNS}$) without transfer learning (top row) and with transfer learning (bottom row) for Trajectory~2 in Table~\ref{tab:single_orbit} with $L=25$ (left column) and $L=28$ (right column).}

\begin{table}[ht]
\centering
\caption{\revision{Transfer learning improves single-orbit prediction for the KS model.
We initially trained the ESN using trajectories from the $L=22$ regime.
We used this ESN to predict $5$ trajectories for each integer value of $L$ from $L=23$ to $L=29$.
The top grid reports Lyapunov times for which relative $L^{2}$-errors for these single-orbit predictions remain below $10\%$ (no transfer learning).
After updating the training of the original ESN at a transfer learning level of $5\%$, ESN predictions typically remain accurate for $1$ to $2$ additional Lyapunov times.
Bold indicates strict increase.}}
\label{tab:single_orbit}
\begin{tabular}{| c  |c  |c  |c  |c  |c  |c  |c |}
\hline
& $L=23$ & $L=24$ & $L=25$ & $L=26$ & $L=27$ & $L=28$ & $L=29$ \\
\hline
Trajectory $1$ & 11& 10& 9& 7& 6& 3& 2\\
\hline
Trajectory $2$ & 10 & 9 & 9 & 8 & 6 & 4 & 2 \\
\hline
Trajectory $3$ & 11 & 10 & 10 & 9 & 8 & 6 & 4 \\
\hline
Trajectory $4$ & 9 & 9 & 8 & 6 & 5 & 3 & 2 \\
\hline
Trajectory $5$ & 12 & 11 & 10 & 9 & 6 & 4 & 3 \\
\hline
\end{tabular}
\\
\vspace{3ex}
\begin{tabular}{| c  |c  |c  |c  |c  |c  |c  |c |}
\hline
& $L=23$ & $L=24$ & $L=25$ & $L=26$ & $L=27$ & $L=28$ & $L=29$ \\
\hline
Trajectory $1$ & 11& 10& 9& \textbf{8}& \textbf{7}& \textbf{5}& \textbf{3}\\
\hline
Trajectory $2$ & 10 & \textbf{10} & \textbf{10} & 8 & \textbf{7} & \textbf{6} & \textbf{3} \\
\hline
Trajectory $3$ & 11 & \textbf{10} & 10 & 9 & 8 & \textbf{7} & 4 \\
\hline
Trajectory $4$ & 9 & 9 & \textbf{9} & \textbf{8} & \textbf{6} & \textbf{4} & \textbf{4} \\
\hline
Trajectory $5$ & 12 & 11 & 10 & \textbf{10} & \textbf{8} & \textbf{6} & \textbf{4} \\
\hline
\end{tabular}
\end{table}

\subsection{Predicting the statistical behavior of the gKS equation}
\label{sec:gKS_prediction}

We consider the gKS equation when the parameter $\gamma$ in the dispersion relation is small.
We assume that for $\gamma$ small, the gKS dynamics admit an attractor that supports an ergodic invariant measure.
We assume that this invariant measure describes the asymptotic distribution of all of the initial data that we select.
We show that for $\gamma$ small and fixed, the ESN accurately predicts the power spectrum of the gKS equation \revision{and the correlation dimension of the attractor} for this fixed value of $\gamma$, once trained on trajectories generated using this fixed value of $\gamma$.
We then show that transfer learning is remarkably effective for small $\gamma$.
Throughout Section~\ref{sec:gKS_prediction}, we set $L=43$.


\subsubsection{No transfer learning}

We performed the following experiment for each of two values of $\gamma$, $\gamma = 0$ and $\gamma = 0.1$.
We began with an untrained ESN.
We generated $30$ gKS trajectories, \revision{each of length $T=10000$ (approximately $1000$ Lyapunov times)}, via direct numerical simulation.
We used $20$ of these to train the ESN and the other $10$ for prediction.

\begin{figure}[ht]
\centering
\centerline{
\includegraphics[width=0.48\linewidth]{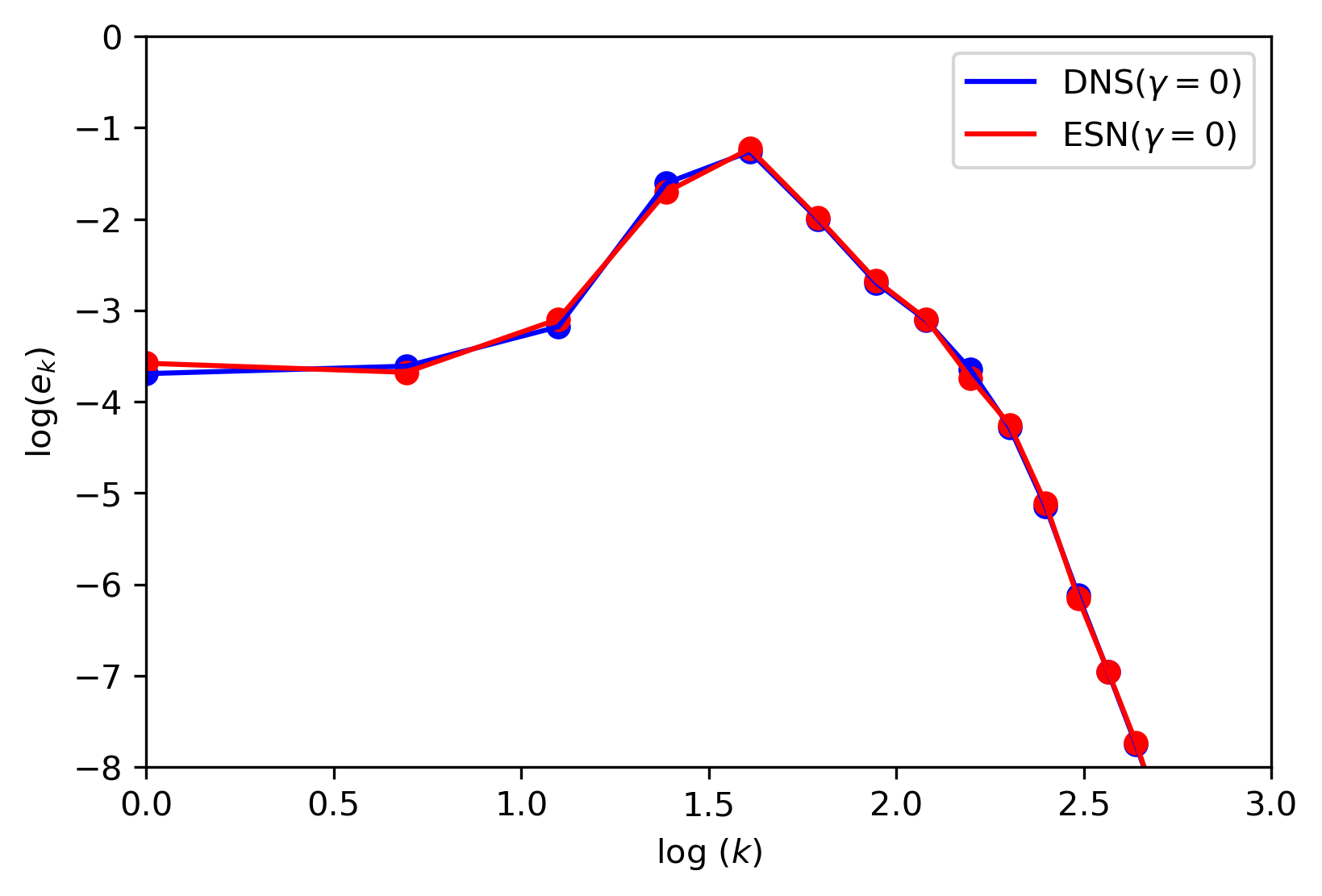} 
\includegraphics[width=0.48\linewidth]{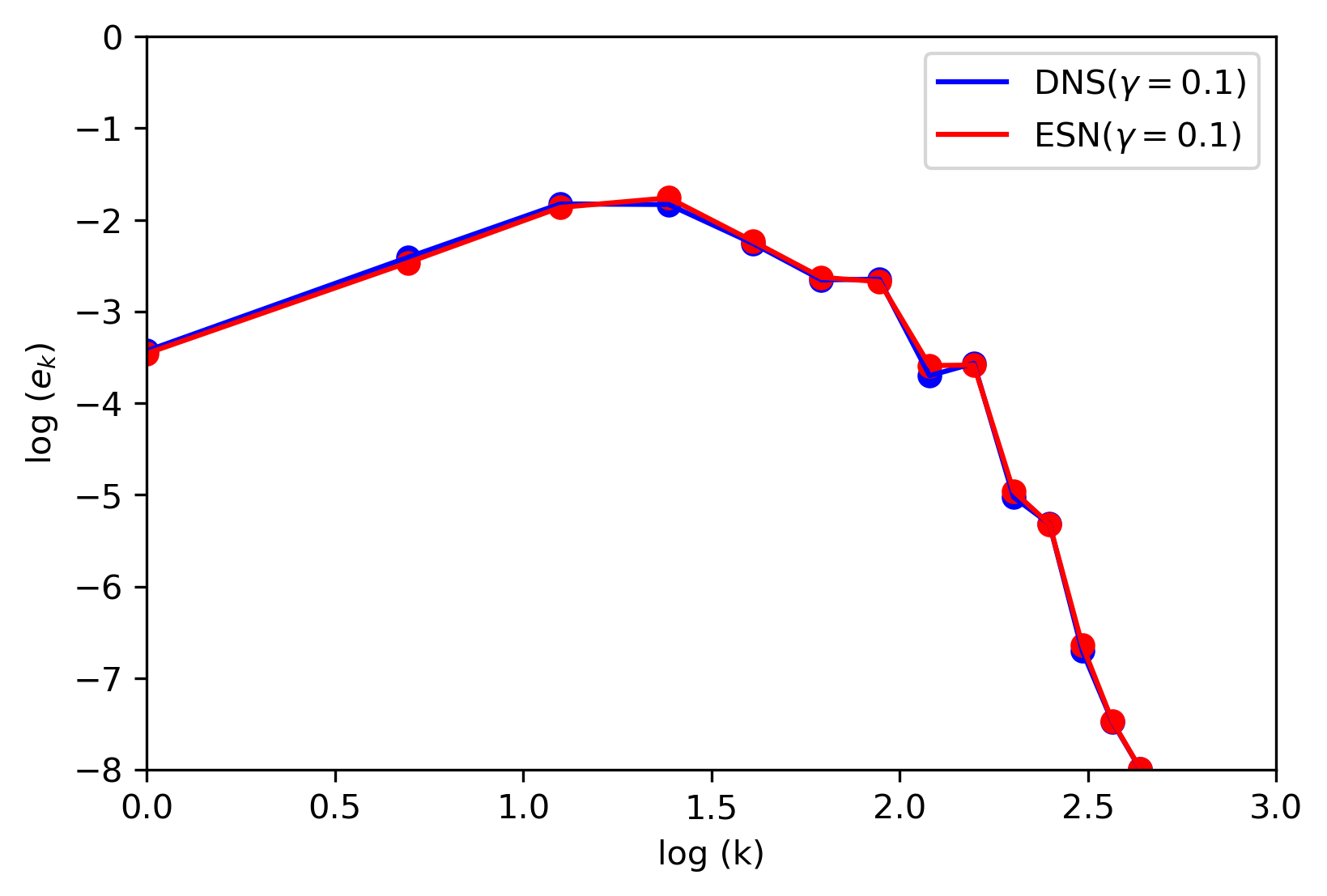}}
\caption{Once trained, the ESN accurately predicts the \revision{\textit{time-averaged}}
power spectrum of the gKS equation.
Left: $\gamma = 0$.
Right: $\gamma = 0.1$.
Each blue curve shows the power spectrum obtained by averaging over the $10$ DNS trajectories (each of length $T=10000$) used for prediction.
Each red curve shows the power spectrum predicted by the ESN, after the untrained ESN has been trained using $20$ trajectories (each of length $T=10000$).
Plots show $\log (e_{k})$ versus $\log (k)$, where $k$ denotes wavenumber.
Spatial domain: $L=43$.}
\label{fig:gamma_no_TL}
\end{figure}

Figure~\ref{fig:gamma_no_TL} illustrates that for $\gamma = 0$ (left) and $\gamma = 0.1$ (right), the trained ESN accurately predicts the power spectrum of the gKS equation.
Each red curve shows the ESN prediction of the power spectrum.
Each blue curve shows the power spectrum obtained by averaging over the $10$ DNS trajectories used for prediction.
\revision{Table~\ref{tab:gKS_no_TL_corrDim} reports that the trained ESN accurately captures the correlation dimension of the attractor as well.}
Importantly, the power spectrum of the gKS equation \revision{and the correlation dimension of the attractor} depend sensitively on $\gamma$ in the chaotic (small $\gamma$) regime.
This observation motivates the use of transfer learning in this regime.

\begin{table}[ht]
\centering
\caption{\revision{Once trained, the ESN accurately predicts the correlation dimension of the attractor of the gKS equation.
These correlation dimensions correspond to Figure~\ref{fig:gamma_no_TL}.}}
\label{tab:gKS_no_TL_corrDim}
\begin{tabular}{| c | c | c |}
\hline
  & \textbf{$\gamma =0$} & \textbf{$\gamma =0.1$}  \\
\hline
DNS & 6.442 & 3.416 \\
\hline
ESN & 6.311 & 3.449 \\
\hline
\end{tabular}
\end{table}


\subsubsection{Using transfer learning for the gKS equation}
\label{sec:gKSTL}

Figure~\ref{fig:gamma_no_TL} demonstrates that the statistical behavior of the gKS equation depends sensitively on $\gamma$ in the small-$\gamma$ regime.
We therefore ask if transfer learning can capture substantial changes in statistical behavior as $\gamma$ varies.
To address this question, we performed the following experiment.
We started with an untrained ESN and first trained it using $20$ trajectories \revision{(each of length $T=10000$)} from the $L=43$ and $\gamma = 0$ regime.
We then updated the training using \revision{some} TL data from the $\gamma = 0.1$ regime at transfer learning rate $\alpha = 0.005$.
\revision{After the training update, we predicted the power spectrum of the gKS equation for $\gamma = 0.1$.}

\begin{figure}[ht]
\centering
\centerline{
\includegraphics[width=0.7\linewidth]{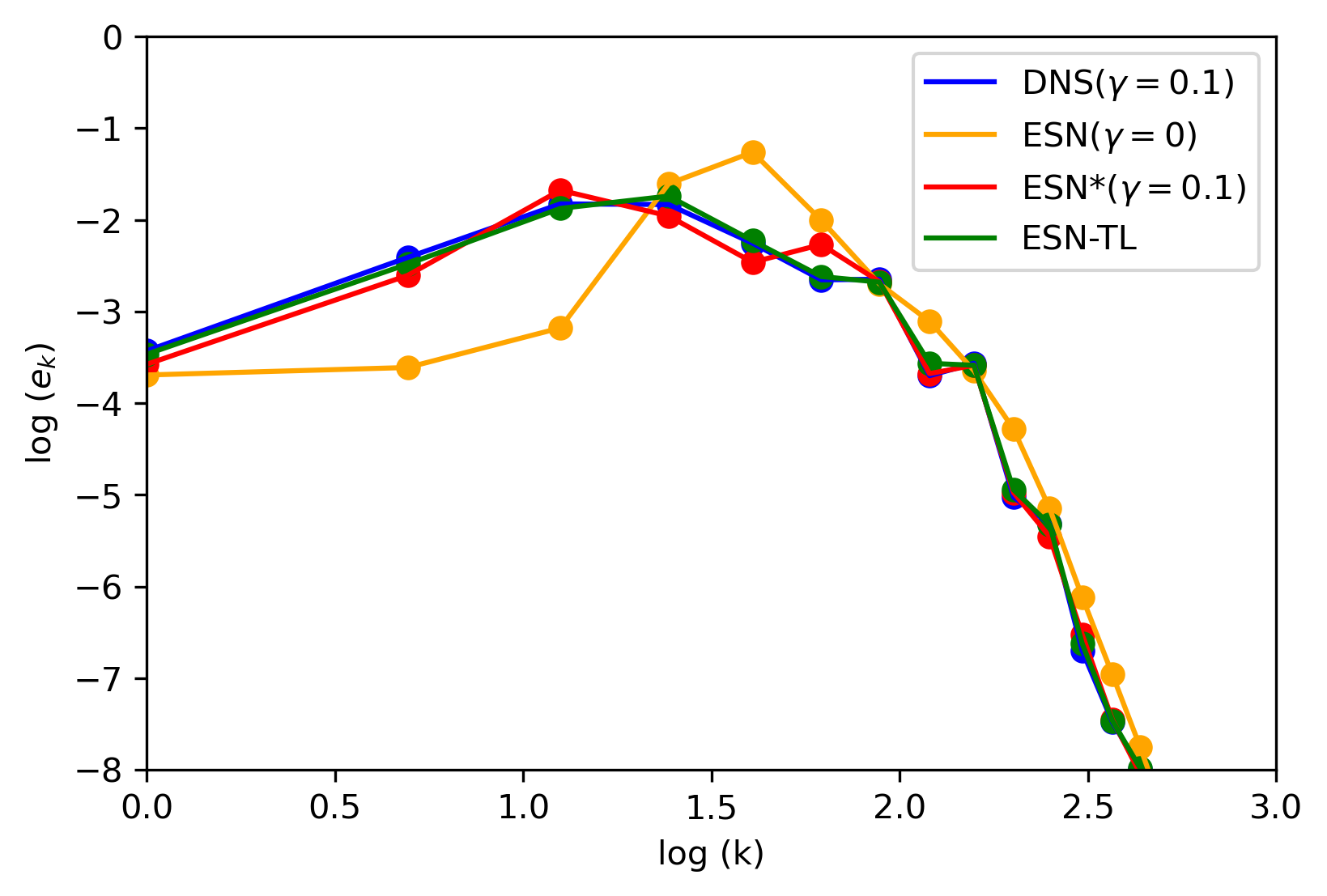}}
\caption{
Transfer learning allows the ESN to accurately predict the \revision{\textit{time-averaged}}
power spectrum of the gKS equation as the dispersion parameter $\gamma$ \revision{instantaneously jumps}.
Blue curve: Power spectrum of the gKS equation in the $\gamma = 0.1$ regime obtained by averaging over the $10$ DNS trajectories \revision{(each of length $T=10000$)} used for prediction.
Green curve: Power spectrum predicted by the ESN after a first round of training using $20$ trajectories from the $\gamma = 0$ regime \revision{(each of length $T=10000$)} and then a training update using TL data from the $\gamma = 0.1$ regime at level $10\%$.
Yellow curve: Power spectrum predicted by the ESN after training only on $20$ trajectories from the $\gamma = 0$ regime \revision{(no transfer learning)}.
Red curve: Power spectrum predicted by an ESN that has been trained only on the TL data, namely $2$ trajectories from the $\gamma = 0.1$ regime.
Plots show $\log (e_{k})$ versus $\log (k)$, where $k$ denotes wavenumber.
Spatial domain: $L=43$.
}
\label{fig:gamma_with_TL}
\end{figure}

\begin{figure}[ht]
\centering
\centerline{
\includegraphics[width=0.7\linewidth]{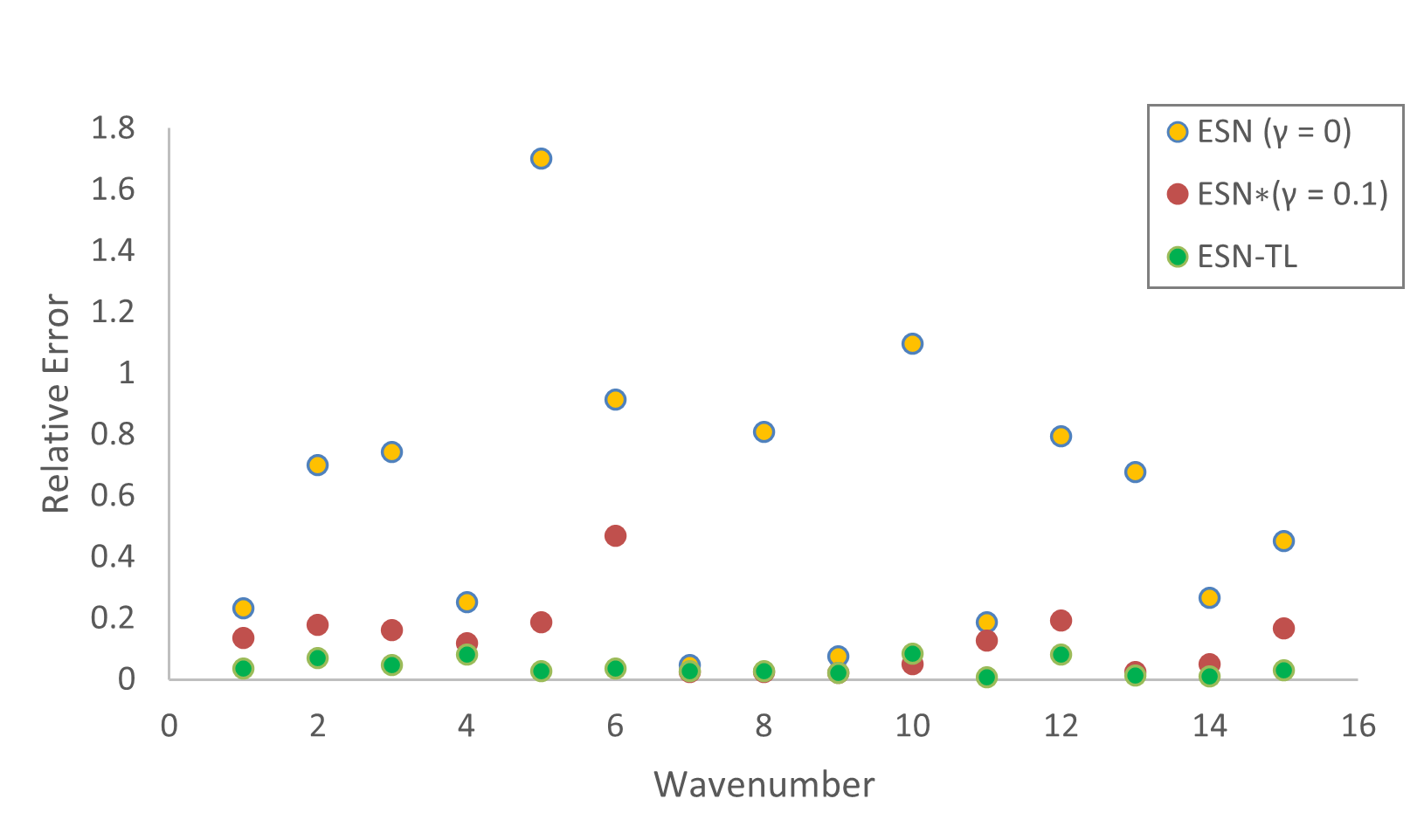}}
\caption{Plot corresponding to Table~\ref{tab:error_gKS_TL}.
}
\label{fig:scatterplot_gamma}
\end{figure}

The result of this experiment is shown in Figure~\ref{fig:gamma_with_TL}.
The blue curve shows the ground truth, namely the power spectrum in the $\gamma = 0.1$ regime obtained by averaging over the $10$ DNS trajectories \revision{(each of length $T=10000$)} from the $\gamma = 0.1$ regime used for prediction.
Using a transfer learning level of only $10\%$ ($2$ trajectories) produces an accurate ESN prediction of the power spectrum (green curve).
For reference, the yellow curve shows the ESN-based prediction of the power spectrum when the ESN is trained only on $20$ trajectories from the $\gamma = 0$ regime \revision{(no transfer learning)}.
\revision{See Table~\ref{tab:error_gKS_TL} \revisiontwo{and Figure~\ref{fig:scatterplot_gamma}} for relative error~\eqref{e:relative_error} as a function of wavenumber.}
\revision{Importantly, the TL-enhanced ESN can continue to accurately predict the power spectrum over long timescales.}
By contrast, single-orbit predictions in the $\gamma = 0.1$ regime based on transfer learning are valid for much shorter timescales~\cite[Section~5.1.1]{shahalam2024}.

As a control, we compare the prediction using transfer learning (Figure~\ref{fig:gamma_with_TL}, green curve, ESN-TL) to a prediction from an ESN that has been trained only on the TL data, namely $2$ trajectories from the $\gamma = 0.1$ regime (Figure~\ref{fig:gamma_with_TL}, red curve, ESN*).
\revision{Visually (Figure~\ref{fig:gamma_with_TL}) and quantitatively (Table~\ref{tab:error_gKS_TL} \revisiontwo{and Figure~\ref{fig:scatterplot_gamma}}), the ESN* prediction is inferior to the ESN-TL prediction.}
Discrepancies between these two predictions occur at low wavenumbers that contain most of the energy.
The curves shown in Figure~\ref{fig:gamma_with_TL} contain total energies $2.92$ (blue, DNS), $2.80$ (red, ESN*), and $2.97$ (green, ESN-TL).
This corresponds to relative error $4.1\%$ and $1.7\%$ for the ESN* and ESN-TL predictions, respectively.
In this sense, the ESN-TL prediction is more than twice as accurate as the ESN* control.

\revision{Table~\ref{tab:gKS_TL_corrDim} reports correlation dimension of the attractor for the experiment from Figure~\ref{fig:gamma_with_TL}.
We see that transfer learning allows the ESN to accurately predict correlation dimension for $\gamma = 0.1$, and the ESN-TL prediction is better than the control (ESN*) prediction.}

\begin{table}[ht]
\caption{\revision{Correlation dimension of the attractor of the gKS equation for the experiment from Figure~\ref{fig:gamma_with_TL}.}}
\label{tab:gKS_TL_corrDim}
\centering
\begin{tabular}{|c|c|c|c|}
\hline
DNS ($\gamma = 0.1$) & ESN ($\gamma = 0$) & ESN* ($\gamma = 0.1$) & ESN-TL \\
\hline
3.416 & 6.311 & 3.242 & 3.401 \\
\hline
\end{tabular}
\end{table}


\section{Discussion}

The interface of machine learning and dynamical systems is a rapidly developing field.
Echo state networks have emerged as a powerful tool.
For spatiotemporally chaotic systems, these networks are capable of both single-orbit and statistical prediction.
In this paper, we have focused on predicting the long-term statistical behavior of an archetypal model, the gKS equation in the spatiotemporally chaotic regime.

We have shown that when the domain size $L$ and the dispersion parameter $\gamma$ are fixed, our ESN can learn and then accurately predict the power spectrum of the gKS equation \revision{and the correlation dimension of the attractor}.
Further, we have shown that transfer learning allows one to efficiently and accurately track changes in the power spectrum of the gKS equation \revision{and the correlation dimension of the attractor} when $L$ or $\gamma$ \revision{instantaneously jumps}.
Importantly, our ESN did not blow up nor collapse during our long simulations.
This is a great sign for the utility of ESNs with respect to long-term statistical prediction.

\revision{Our results suggest that transfer learning is more efficient when the correlation dimension of the attractor drops upon a parameter set jump than when it rises upon such a jump.
In Section~\ref{sec:KSTL}, we initially train the ESN using KS trajectories for $L=22$ and then apply transfer learning to do statistical prediction for the KS equation when $L=43$.
The jump from $L=22$ to $L=43$ induces a jump in the correlation dimension of the attractor (Table~\ref{tab:L_TL_corrDim}).
In this scenario, we observe that a transfer learning level of \revisiontwo{$50\%$} produces accurate statistical predictions.
By contrast, in Section~\ref{sec:gKSTL} we fix $L=43$, initially train when $\gamma = 0$, and then apply transfer learning to do statistical prediction for the gKS equation when $\gamma = 0.1$.
This time, the jump from $\gamma = 0$ to $\gamma = 0.1$ induces a drop in the correlation dimension of the attractor (Table~\ref{tab:gKS_TL_corrDim}).
In this second scenario, we observe that a transfer learning level of only $10\%$ suffices for accurate statistical predictions.
It is intuitive that transfer learning is more efficient when the correlation dimension of the attractor decreases upon a parameter set jump than when it rises upon such a jump, since the output layer $W_{\mathrm{out}}$ must be updated to project higher-dimensional chaos in the latter case.}

Looking to the future, ESNs can potentially be used to inexpensively simulate climate \revision{or to generate chaotic timeseries with known statistical properties.}
Importantly, ESNs allow relatively large timesteps when generating new trajectories, making \revision{ESNs} more efficient than traditional numerical integrators.
Thus far, transfer learning for ESNs has received limited attention in the literature.
We believe transfer learning is an important concept that can be applied to complex problems and allow ESNs to track parametric changes in spatiotemporally chaotic dynamics.
In general, developing efficient models that can track parametric changes in dynamical systems is an important task.
\revision{The methodology described in this paper aims at developing efficient computational ML models for this task.}

\section{Code}

Our code for simulating the gKS equation and ESNs is located on GitHub~\cite{shah2025esn}.

\section{Author contributions}

All authors wrote the paper.
Study design: WO, IT.
Initial code development: IT.
Code development, simulation, and visualization: MSA.

\section{Supporting tables and figures}

\begin{table}[ht]
\centering
\caption{\revision{Relative error~\eqref{e:relative_error} as a function of wavenumber for the TL-based predictions in Figure~\ref{fig:L_with_TL}.
}}
\label{tab:error_L_TL}
\begin{tabular}{| c | c | c | c | c |}
\hline
Wavenumber & ESN-L22 & ESN-TL 10\% & ESN-TL 25\% & ESN-TL 50\% \\
\hline
1 & 0.603 & 0.445 & 0.194 & 0.048 \\
\hline
2 & 10.496 & 0.771 & 0.307 & 0.131 \\
\hline
3 & 4.403 & 1.022 & 0.248 & 0.014 \\
\hline
4 & 0.468 & 0.336 & 0.101 & 0.034 \\
\hline
5 & 0.943 & 0.787 & 0.381 & 0.004 \\
\hline
6 & 0.965 & 0.495 & 0.118 & 0.004 \\
\hline
7 & 0.989 & 0.536 & 0.211 & 0.128 \\
\hline
8 & 0.997 & 0.319 & 0.072 & 0.113 \\
\hline
9 & 0.999 & 0.308 & 0.237 & 0.047 \\
\hline
10 & 1.000 & 0.465 & 0.137 & 0.026 \\
\hline
11 & 1.000 & 0.416 & 0.031 & 0.013 \\
\hline
12 & 1.000 & 0.303 & 0.197 & 0.054 \\
\hline
13 & 1.000 & 0.324 & 0.111 & 0.062 \\
\hline
14 & 1.000 & 0.298 & 0.064 & 0.072 \\
\hline
15 & 1.000 & 0.585 & 0.132 & 0.035 \\
\hline
\end{tabular}
\end{table}

\begin{table}[ht]
\centering
\caption{\revision{Relative error~\eqref{e:relative_error} as a function of wavenumber for the ergodicity experiment described in Section~\ref{sec:KSTL}.
}}
\label{tab:error_ergodicity}

\begin{tabular}{| c | c | c | c | c |}
\hline
Wavenumber & ESN-L22 & ESN-TL 10\% & ESN-TL 25\% & ESN-TL 50\% \\
\hline
1 & 0.611& 0.432 & 0.193 & 0.046 \\
\hline
2 & 10.315& 0.770 & 0.313 & 0.136 \\
\hline
3 & 4.372& 1.005 & 0.249 & 0.014 \\
\hline
4 & 0.459& 0.348 & 0.102 & 0.033 \\
\hline
5 & 0.927& 0.766 & 0.377 & 0.005 \\
\hline
6 & 0.971& 0.480 & 0.122 & 0.013 \\
\hline
7 & 0.982& 0.518 & 0.211 & 0.125 \\
\hline
8 & 0.998& 0.316 & 0.073 & 0.112 \\
\hline
9 & 1.000& 0.316 & 0.242 & 0.047 \\
\hline
10 & 1.000 & 0.463 & 0.140 & 0.025 \\
\hline
11 & 1.000 & 0.429 & 0.031 & 0.013 \\
\hline
12 & 1.000 & 0.294 & 0.198 & 0.055 \\
\hline
13 & 1.000 & 0.325 & 0.110 & 0.062 \\
\hline
14 & 1.000 & 0.305 & 0.064 & 0.071 \\
\hline
15 & 1.000 & 0.590 & 0.135 & 0.036 \\
\hline

\end{tabular}

\end{table}

\begin{table}[ht]
\centering
\caption{\revision{Relative error~\eqref{e:relative_error} as a function of wavenumber for the experiment in Figure~\ref{fig:gamma_with_TL}.
}}
\label{tab:error_gKS_TL}
\begin{tabular}{| c | c | c | c |}
\hline
Wavenumber & ESN $(\gamma =0)$ & ESN$^*(\gamma=0.1)$ & ESN-TL \\
\hline
1 & 0.232& 0.137& 0.036\\
\hline
2 & 0.699& 0.178& 0.07\\
\hline
3 & 0.742& 0.161& 0.048\\
\hline
4 & 0.253& 0.118& 0.081\\
\hline
5 & 1.702& 0.187& 0.028\\
\hline
6 & 0.915& 0.470& 0.036\\
\hline
7 & 0.048& 0.025& 0.028\\
\hline
8 & 0.808& 0.023& 0.027\\
\hline
9 & 0.077& 0.021& 0.021\\
\hline
10 & 1.097& 0.050& 0.083\\
\hline
11 & 0.188& 0.128& 0.006\\
\hline
12 & 0.795& 0.192& 0.08\\
\hline
13 & 0.676& 0.025& 0.012\\
\hline
14 & 0.267& 0.051& 0.009\\
\hline
15 & 0.452& 0.167& 0.031\\
\hline
\end{tabular}
\end{table}

\FloatBarrier


\end{document}